\documentclass[11pt]{amsart}

\usepackage{amsfonts, amssymb, amscd}
\numberwithin{equation}{section}

\usepackage[symbol]{footmisc}

\usepackage{bm}
\usepackage{verbatim}
\usepackage{amssymb}
\usepackage{mathrsfs}
\usepackage{graphicx}
\usepackage{tikz-cd}
\usepackage{subcaption}
\usepackage{subfiles}
\usepackage[toc,page]{appendix}
\usepackage{mathtools}
\usepackage{comment}
\usepackage{enumerate}
\usepackage{enumitem}

\usepackage{graphicx}
\graphicspath{{images/}}

\usepackage{appendix}
\usepackage{hyperref}

\newcommand{\id}{\mathrm{id}}

\newcommand{\Qq}{\mathbb{Q}}

\newcommand{\vol}{\operatorname{vol}}

\newcommand{\rk}{\operatorname{rank}}

\newcommand{\tang}{\operatorname{tang}}
\newcommand{\Res}{\operatorname{Res}}
\newcommand{\Null}{\operatorname{Null}}

\newcommand{\Ex}{\operatorname{Ex}}
\newcommand{\Supp}{\operatorname{Supp}}

\newcommand{\Bs}{\operatorname{Bs}}
\newcommand{\const}{\operatorname{const}}
\newcommand{\mult}{\operatorname{mult}}

\newcommand{\Ff}{\mathscr{F}}

\newcommand{\Sing}{\mathrm{Sing}}

\newtheorem{thm}{Theorem}[section]

\newtheorem{cor}[thm]{Corollary}
\newtheorem{lem}[thm]{Lemma}
\newtheorem{defn}[thm]{Definition}

\newtheorem{claim}[thm]{Claim}
\theoremstyle{definition}
\newtheorem{rem}[thm]{Remark}

\theoremstyle{definition}

\begin{document}

\title{On invariance of plurigenera for foliated surface pairs}

\author{Jihao Liu}

\address{Department of Mathematics, The University of Uath, Salt Lake City, UT 84112, USA}
\email{jliu@math.utah.edu}

\begin{abstract} We prove an invariance of plurigenera for some foliated surface pairs of general type.
\end{abstract}

\date{\today}
\maketitle
\pagestyle{myheadings}\markboth{\hfill  Jihao Liu \hfill}{\hfill On invariance of plurigenera for foliated surface pairs\hfill}

\tableofcontents
\section{Introduction}We work over the field of complex numbers $\mathbb C$.

In \cite{CF15}, Cascini and Floris prove an invariance of plurigenera for a family of foliated surfaces:

\begin{thm}[Invariance of plurigenera for foliated surfaces, {\cite[Corollary \ref{thm: foliated surface pair inv plur}]{CF15}}]\label{thm: foliated surface inv plur}  Let $(X_t, \Ff_t)|_{t\in T}$ be a smooth family of foliations of surfaces with at most reduced singularities. Then for any sufficiently large integer $m>0$, 
$$h^0(X_t, mK_{\Ff_t})=\const$$ 
for any $t\in T$, i.e. the invariance of plurigenera $h^0(X_t, mK_{\Ff_t})$ holds for sufficiently large $m$.
\end{thm}

In birational geometry, it is usually natural and useful to consider pairs $(X,\Delta)$, where $X$ is a normal variety and $\Delta\geq 0$ is an effective $\Qq$-divisor such that $K_X+\Delta$ is $\Qq$-Cartier. Therefore, for normal varieties associated with foliations, we may similarly consider \emph{foliated pairs}:
\begin{defn}[See Definition \ref{defn: foliated triples}]\label{defn: foliated triples introduction}
A \emph{foliated triple} $(X,\Ff,\Delta)$ consists of a normal variety $X$, a foliation $\Ff$ on $X$ an an effective $\Qq$-divisor $\Delta$ on $X$ such that $K_{\Ff}+\Delta$ is $\Qq$-Cartier. 

We usually say that $(\Ff,\Delta)$ is a \emph{foliated pair} on $X$. 
\end{defn}

It is then natural to ask whether Theorem \ref{thm: foliated surface inv plur} holds for a family of foliated surface pairs. Part of the intuition of this question comes from the invariance of plurigenera for klt pairs, which is an important generalization of the well-known Siu's deformation invariance of plurigenera (cf.~\cite{Siu98}).
\begin{thm}[Deformation invariance of log plurigenera, {\cite[Theorem 1.8(1)]{HMX13}}]\label{thm: deformation invariance plurigenera klt pair} 
Let $\pi: X\rightarrow T$ be a projective morphism of smooth varieties. Suppose that $(X,\Delta)$ is klt with simple normal crossings over $T$ and either $\Delta$ or $K_X+\Delta$ is big. Then $h^0(X_t,m(K_{X_t}+\Delta_t))$ is independent of $t\in T$ for any integer $m>0$ such that $m\Delta$ is integral.
\end{thm}

In this paper we prove the following result:
\begin{thm}\label{thm: foliated surface pair inv plur} Let $(X_t, \Ff_t)|_{t\in T}$ be a smooth family of canonical foliations of surfaces and $\Delta$ an effective $\mathbb Q$-Cartier $\mathbb Q$-divisor on $X$, such that 
\begin{itemize}
\item $K_\Ff+\Delta$ is big over $T$, and
\item each irreducible component of $\Delta$ is irreducible and reduced over any closed point of $T$, i.e. if we write $\Delta=\sum a_i\Delta_i$ into its irreducible components, then each $\Delta_i|_{X_t}$ is irreducible and reduced for any closed point $t\in T$,
\end{itemize}
then for any sufficiently divisible integer $m>0$,
$$h^0(X_t, m(K_{\Ff_t}+\Delta_t))=\const.$$
\end{thm}

As a corollary, we deduce the invariance of volumes of foliated surface pairs for a smooth family of canonical foliations:
\begin{cor}\label{cor: foliated surface pair vol inv} Let $(X_t, \Ff_t)|_{t\in T}$ be a smooth family of canonical foliations of surfaces. For any effective $\mathbb Q$-Cartier $\mathbb Q$-divisor $\Delta$ on $X$  such that each irreducible component of $\Delta$ is irreducible and reduced over any closed point of $T$, we have
$$\vol(X_t,K_{\Ff_t}+\Delta_t)=\const.$$
\end{cor}

A by-product in the proof of Theorem \ref{thm: deformation invariance plurigenera klt pair} is a vanishing theorem for foliated surface pairs. 

\begin{thm}\label{thm: weak kv vanishing foliated surface} Let $X$ be a smooth surface, $\Ff$ a rank 1 canonical foliation on $X$, and $\Delta$ an effective $\mathbb Q$-Cartier $\mathbb Q$-divisor on $X$, such that 
\begin{itemize}
\item $K_\Ff+\Delta$ is big and nef,
\item $\lfloor\Delta\rfloor=0$, and
\item each irreducible component of $\Delta$ is not $\Ff$-invariant.
\end{itemize}
Then 
$$H^1(X,m(K_\Ff+\Delta))=0$$ 
for any sufficiently divisible integer $m>0$.
\end{thm}

\noindent\textit{Sketch of the paper}. In Section 2, we introduce basic notation and conventions of foliations. In Section 3 we review the $K_\Ff$-MMP for foliated surfaces. In Section 4, we give a precise description on the $(K_\Ff+\Delta)$-MMP for certain foliated pairs. In Section 5, we study the structure of the Zariski decomposition for certain foliated surface pairs, and prove a vanishing theorem (Lemma \ref{lem: lem: kv vanishing foliation version 2}). We prove Theorem \ref{thm: weak kv vanishing foliated surface} at the end of this section. In Section 6, we study the deformation property of foliated pairs as in \cite[Chapter 4]{Bru01}, and show that the negative part of the Zariski decomposition of $K_\Ff+\Delta$ deforms. In Section 7, we use technical strategies similar to \cite[Section 4]{HMX13} and the vanishing theorem proved in Section 5 to prove our main theorems.\vspace{2mm}

%


%

\noindent\textbf{Acknowledgement}. The author would like to thank Christopher D. Hacon for his constant support, encouragement and inspiring discussions through the writing of the paper. The author also wishes to thank Enrica Floris, Jorge V. Pereira, and Chenyang Xu for their carefully reading the manuscripts and giving numerous helpful comments. The author was supported by NSF research grants no: DMS-1300750, DMS-1265285 and by a grant from the Simons Foundation (Award number: 256202).

\section{Notation and conventions}

In this section, we give the definitions which we frequently use in this paper. Most of these definitions are adopted from \cite{Bru00},~\cite{Bru01} and \cite{Bru02}, but the readers should be careful as there may exist small differences.

\subsection{Foliations}

\begin{defn}[Foliation]\label{defn: foliation} Let $X$ be a normal variety.  A \emph{foliation} $\Ff$ on $X$ is a coherent saturated subsheaf of the tangent sheaf $T_X$ of $X$ such that $\Ff$ is closed under the Lie bracket. 

If $X$ is associated with a foliation $\Ff$ on $X$, then $X$ is called a \emph{foliated variety}. A \emph{foliated surface} is a foliated variety of dimension $2$. 
\end{defn}

\begin{defn}[Singularities of foliation]\label{defn: sing foliation}
Let $X$ be a normal variety and $\Ff$ a foliation on $X$. The \emph{singular locus} of $\Ff$ is defined as the locus where $\Ff$ fails to be a sub-bundle of $T_{X}$, and is denoted by $\Sing(\Ff)$. 

If $\Sing(\Ff)=\emptyset$, $\Ff$ is called \emph{smooth}. For any closed point $p\in X$, $p$ is called a \emph{singular point of $\Ff$} if $p\in\Sing(\Ff)$, and $p$ is called a \emph{smooth point of $\Ff$} if $p\not\in\Sing(\Ff)$.

Let $U=X\backslash \Sing(X)$ be the smooth locus of $X$. The \emph{canonical divisor} of $\Ff$ is the divisor associated with $\det(\Ff|_U)^*$, and is denoted by $K_\Ff$. 
\end{defn}

\begin{rem}[Abuse of notation]\label{rem: sing foliation not sing X}
In the rest of the paper, we only consider foliations which additionally satisfy the following property: suppose that $\Ff$ is a foliation on a normal variety $X$, then
\begin{equation}\label{equ: sing not intersect}
\Sing(\Ff)\cap\Sing(X)=\emptyset.\tag{$\star$}
\end{equation}
Although many things can be done in greater generality, this assumption follows the terminology as in \cite{Bru02} in order to define algebraic invariants for foliated surfaces with at most cyclic quotient singularities. Nevertheless, this property always holds when $X$ is smooth.
\end{rem}

\begin{rem}\label{rem: sing foliation codim 2}
By definition, $\Sing(\Ff)$ is a codimension $2$ subset of $X$. In particular, if $\dim X=2$, $\Sing(\Ff)$ is a set of isolated closed points. We will repeatedly use this fact in the rest of this paper.
\end{rem}

\begin{defn}[Restriction, pullback and pushforward of foliations]\label{defn: basic property foliation} Let $X$ be a normal variety and $\Ff$ a foliation on $X$. We give the following definitions.
\begin{itemize}
\item For any open subset $U\subset X$, the \emph{restriction} of $\Ff$ to $U$ is defined as the foliation associated with the restriction sheaf $\Ff|_U$, and is also denoted by $\Ff|_U$.
\item For any birational morphism $\pi: Y\rightarrow X$, $\Ff$ induces a foliation on $Y\backslash\Ex(f)$, which can be extended to a foliation $\Ff_Y$ on $Y$. $\Ff_Y$ is called the \emph{pullback} of $\Ff$ on $Y$ and is denoted by $\pi^*\Ff$.
\item For any birational morphism $\nu: X\rightarrow Z$ such that the sheaf $\nu_*\Ff$ is a foliation on $Z$, the \emph{pushforward} of $\Ff$ to $Z$ is the foliation associated with $\nu_*\Ff$, and is also denoted by $\nu_*\Ff$.
\end{itemize}
\end{defn}

\begin{defn}[Canonical foliation]\label{defn: canonical sing foliation} Let $X$ be a normal variety and $\Ff$ a foliation on $X$. 
\begin{itemize}
    \item $\Ff$ is called \emph{canonical}, if for any birational morphism $\pi: Y\rightarrow X$, 
$$K_{\pi^*{\Ff}}=\pi^*{K_{\Ff}}+E,$$
for some effective $\Qq$-divisor $E$ on $Y$.
\item A \emph{canonical singularity} of $\Ff$ is a closed point $p\in\Sing(\Ff)$, such that there exists an open neighborhood $U$ of $p$ such that $\Ff|_{U}$ is canonical.
\end{itemize}
\end{defn}

\begin{defn}[Rank of a foliation]\label{defn: rk foliation}
Let $X$ be a normal variety and $\Ff$ a foliation on $X$. The \emph{rank} of $\Ff$ is defined as the rank of $\Ff$ as a sheaf, and is denoted by $\rk(\Ff)$. The \emph{co-rank} of $\Ff$ is defined as $\dim X-\rk(\Ff)$. 
\end{defn}

\begin{defn}[Rank 1 foliations and co-rank 1 foliations]\label{defn: rk 1 foliation}
Let $X$ be a normal variety, $\Ff$ a foliation on $X$ and $p\in X\backslash \Sing(X)$ a closed point. It is immediate from the definitions that we have the following:
\begin{itemize}
\item If $\Ff$ is a co-rank 1 foliation, then there exists an open neighborhood $U$ of $p$ such that $\Ff$ is defined by the kernel of a holomorphic 1-form $\omega$ on $U$. For every such $\omega$, we say that $\omega$ \emph{generates} $\Ff$ around $p$.
\item If $\Ff$ is a rank 1 foliation, then there exists an open neighborhood $U$ of $p$, such that $\Ff$ is defined by a holomorphic vector field $v$ on $U$. For every such $v$, we say that $v$ \emph{generates} $\Ff$ around $p$.
\item In particular, if $X$ is a surface and $\Ff$ is a rank 1 foliation, then $\Ff$ is both generated by a holomorphic-1 form $\omega$ around $p$ and generated by the vector field $v:=\ker(\omega)$ around $p$.
\end{itemize}
\end{defn}

The next easy lemma shows that the smooth points of a foliation are canonical:
\begin{lem}[Co-rank 1 smooth foliation is canonical,{~\cite[Lemma 2.4]{Spi17}}]\label{lem: co-rank 1 smooth is canonical}
Let $X$ be a normal variety, $\Ff$ a co-rank 1 foliation on $X$, and $x\notin \Sing(X)$ a smooth point of $\Ff$. Then $\Ff$ is canonical near an open neighborhood of $x$.
\end{lem}

At the end of this subsection, we give the definitions of \emph{foliated triples} and \emph{foliated pairs}.
\begin{defn}\label{defn: foliated triples}
A \emph{foliated triple} $(X,\Ff,\Delta)$ consists of a normal variety $X$, a foliation $\Ff$ on $X$ and an effective $\Qq$-divisor $\Delta$ on $X$ such that $K_{\Ff}+\Delta$ is $\Qq$-Cartier. We usually call $(\Ff,\Delta)$ a \emph{foliated pair} on $X$. 

A \emph{morphism} between two foliated triples $(X,\Ff,\Delta)$ and $(X',\Ff',\Delta')$ is a morphism $\phi: X\rightarrow X'$, such that $\phi_*\Ff=\Ff'$ and $\phi_*\Delta=\Delta'$. We usually use the notation $\phi: (X,\Ff,\Delta)\rightarrow (X',\Ff',\Delta')$.
\end{defn}

\subsection{Foliations on surfaces} 
\begin{rem}[Abuse of notation]\label{rem: foliation surface is rank 1}
In the rest of the paper, if there is no confusion, when we consider foliations on surfaces, we only consider the foliations of rank $1$. 
\end{rem}

\begin{defn}[Surface with at most cyclic quotient singularities]\label{defn: cyc sing foliation} Let $X$ be a normal surface. We say that $X$ \emph{has at most cyclic quotient singularities} if for any $p\in\Sing(X)$, $p$ is of type $\textbf{\rm \textbf{B}}^2/\Gamma_{k,h}$ for some coprime positive integers $k$ and $h$ such that $kh\not=1$. Here $\textbf{\rm \textbf{B}}^2$ means the unit ball in $\mathbb{C}^2$, and $\Gamma_{k,h}$ is the cyclic group acting on $\mathbb{C}^2$ generated by 
$$(z,w)\rightarrow (e^{\frac{2\pi i}{k}}z,e^{\frac{2h\pi i}{k}}w).$$
For simplicity, we define $\mathcal{S}_{cyc}$ to be the set of all the normal projective surfaces with at most cyclic quotient singularities.
\end{defn}

Surface with at most cyclic quotient singularities are always $\Qq$-factorial:
\begin{lem}[{\cite[Chapter 1, p3]{Bru02}}]\label{lem: qfactoriality cyclic quotient surface}
Let $X\in\mathcal{S}_{cyc}$ be a surface and $\Ff$ a foliation on $X$. Then $X$ is $\mathbb Q$-factorial, hence for any foliation $\Ff$ on $X$,  $K_\Ff$ is $\mathbb{Q}$-Cartier. Moreover, if $X$ is smooth, $K_\Ff$ is Cartier.
\end{lem}

\begin{defn}[Eigenvalues, reduced singularities and Poincar\'{e}-Dulac singularities]\label{defn: eigenvalue of foliation} Let $X$ be a projective surface, $\Ff$ a foliation on $X$, $p\in \Sing(\Ff)$ a closed point, and $v$ a vector field on $X$ which generates $\Ff$ around $p$. By \cite[Chapter 1, p10]{Bru00}, $(Dv)|_p$ has two eigenvalues $\lambda_1,\lambda_2$ which do not depend on the choice of $v$. 

$\lambda_1$ and $\lambda_2$ are called the \emph{eigenvalues} of $\Ff$ at $p$. The \emph{quotient of eigenvalues} of $\Ff$ at $p$ is defined as 
$$\lambda=\lambda(\Ff,p):=\frac{\lambda_1}{\lambda_2}\in\mathbb C\cup\{\infty\}$$ 
under the equivalence $0\sim\infty$ and $x\sim\frac{1}{x}$ for any $x\in\mathbb C$.


If $\lambda\not\in\mathbb Q^+$, $p$ is called a \emph{reduced} singularity of $\Ff$. If $\lambda\in \mathbb{N}^+\cup\frac{1}{\mathbb{N}^+}$, $p$ is called a \emph{Poincar\'{e}-Dulac} singularity of $\Ff$. We say that $\Ff$ has \emph{at most reduced singularities} if any $p\in\Sing(\Ff)$ is a reduced singularity of $\Ff$.
\end{defn}

\begin{lem}[{\cite[Chapter 1, p16]{Bru00}}]\label{lem: type of pd sing}
Let $X$ be a normal surface, $\Ff$ a foliation on $X$ and $p\in\Sing(\Ff)$ a Poincar\'{e}-Dulac singularity of $\Ff$. Then $p$ is of one of the following two types (i.e.  one of the following two holomorphic 1-forms with local coordinates $w$ and $z$ generates $\Ff$ around $p$).
\begin{enumerate}
\item[\textbf{\rm \textbf{Type A}}.] $\lambda wdz-zdw$.
\item[\textbf{\rm \textbf{Type B}}.] $(\lambda w+z^\lambda)dz-zdw$.
\end{enumerate}
\end{lem}

\begin{lem}[{\cite[Chapter 8, p114-115]{Bru00}},{~\cite[Proposition 3.3]{LPT11}}]\label{lem: either type b pd or reduced} Let $X$ be a normal surface, $\Ff$ a foliation on $X$, and $p\in\Sing(\Ff)$ a closed point. Then $p$ is a canonical singularity of $\Ff$ if and only if $p$ is a reduced singularity or a \textbf{\rm \textbf{Type B}} Poincar\'{e}-Dulac singularity of $\Ff$. 
\end{lem}

\begin{defn}[$\Ff$-invariant curve]\label{defn: F invariant curve} Let $X$ be a normal projective surface, $\Ff$  a foliation on $X$, and $U:=X\backslash \Sing(X)$ the smooth locus of $X$. A curve $C$ on $X$ is called $\Ff$-\emph{invariant}, if the inclusion map 
$${T_{\Ff|_{U}}}|_C\rightarrow T_{U}|_C$$ 
factors through $T_{C|_{U}}$.\end{defn}

\subsection{Invariants of foliations on surfaces} 

\begin{defn}\label{defn: tang of foliation} Let $X\in\mathcal{S}_{cyc}$ be a surface and $C$ a curve on $X$ such that each irreducible component of $C$ is \textbf{not} $\Ff$-invariant. For any closed point $p\in C$, we define $\tang(\Ff,C,p)$ in the following way. 
\begin{itemize}
\item If $p\notin \Sing(X)$, let $v$ be a vector field generating $\Ff$ around $p$, and $f$ a holomorphic function generating $C$ around $p$. We define 
$$\tang(\Ff,C,p):=\dim_{\mathbb{C}}\frac{\mathscr{O}_p}{<f, v(f)>}.$$
\item If $p\in\Sing(X)$,  suppose that $p$ is a cyclic quotient singularity of order $k$ of type $\textbf{\rm \textbf{B}}^2/\Gamma$. Let $U$ be an open neighborhood of $p$ and $\rho: \textbf{\rm \textbf{B}}^2\rightarrow\textbf{\rm \textbf{B}}^2/\Gamma\cong U$ the cyclic cover around $p$. Let $\tilde p:=\rho^{-1}(p)$, $\widetilde C:=\rho^*(C|_U)$, and $\widetilde{\Ff}$ the foliation induced by the sheaf $\rho^*(\Ff|_U)$. Then $\tilde p$ is a smooth point of $\widetilde{\Ff}$, and we define
$$\tang(\Ff,C,p):=\frac{1}{k}\tang(\widetilde{\Ff},\widetilde C,\tilde p).$$
\end{itemize}

By \cite[Chapter 2.2, p22]{Bru00}, $\tang(\Ff,C,p)\geq 0$ and $\tang(\Ff,C,p)$ is independent of the choice of $v$ and $f$. Moreover, if $\Ff$ is transverse to $C$ at $p\in C$, then $\tang(\Ff,C,p)=0$.  Therefore, the sum
$$\sum_{p\in C}\tang(\Ff,C,p)$$ 
is a finite sum. We define $\tang(\Ff,C):=\sum_{p\in C}\tang(\Ff,C,p).$
\end{defn}

\begin{defn}\label{defn: z and cs of foliation} Let $X\in\mathcal{S}_{cyc}$ be a surface and $C$ a curve on $X$ such that each irreducible component of $C$ is $\Ff$-invariant. For any closed point $p\in C$, we define $Z(\Ff,C,p)$ and $CS(\Ff,C,p)$ in the following way.
\begin{itemize}
\item If $p\notin \Sing(X)$, let $\omega$ be a $1$-form generating $\Ff$ around $p$, and $f$ a holomorphic function generating $C$ around $p$. Then there are uniquely determined holomorphic functions $g,h$ and a holomorphic $1$-form $\eta$ on $X$, such that
\begin{itemize}
    \item $g\omega=hdf+f\eta$, and
    \item $f,h$ are coprime.
\end{itemize}
We define 
$$Z(\Ff,C,p):=\text{ the vanishing order of } \frac{h}{g}|_C \text{ at } p,$$
and
$$CS(\Ff,C,p):=\text{ residue of }-\frac{1}{h}\eta|_C \text{ at } p.$$
By \cite[Chapter 2, p24]{Bru00}, \cite[Chapter 3, p26]{Bru00} and \cite[Chapter 2]{Bru02}, $Z(\Ff,C,p)$ and $CS(\Ff,C,p)$ are independent of the choice of $\omega$.
\item If $p\in C\cap \Sing(X)$, we define 
$$Z(\Ff,C,p)=CS(\Ff,C,p):=0.$$
\end{itemize}
\end{defn}

\begin{lem}\label{lem: basic properties of z and cs} Let $X$ be a normal projective surface, $\Ff$ a foliation on $X$ and $C$ an irreducible curve on $X$. Then for any closed point $p\in C$ such that $p\notin \Sing(X)$,
\begin{enumerate}
\item $Z(\Ff,C,p)$ is an integer,
\item if $C$ is smooth at $p$, then $Z(\Ff,C,p)\geq 0$,
\item if $p\notin \Sing(\Ff)$, then $Z(\Ff,C,p)=CS(\Ff,C,p)=0$, and
\item if $p$ is a canonical singularity of $\Ff$, then 
\begin{itemize}
\item $Z(\Ff,C,p)\geq 0$, 
\item  if $C$ is smooth at $p$, then $Z(\Ff,C,p)\geq 1$, and 
\item  if $C$ is smooth at $p$ and $\lambda(\Ff,p)\not=0$, then $Z(\Ff,C,p)=1$.
\end{itemize}
\end{enumerate}

\end{lem}

\begin{proof} Because $Z(\Ff,C,p)$ is the vanishing order of a rational function, it must be an integer, hence (1). If $C$ is smooth at $p$, let $v$ be a local holomorphic vector field generating $\Ff$ around $p$. Then $Z(\Ff,C,p)$ is the vanishing order of $v|_C$ at $p$, which is non-negative, hence (2).

Since $C$ is an $\Ff$-invariant curve, $T_{\Ff|_{U}}|_C\rightarrow T_{U}|_C$ factors through $T_{C|_{U}}$. Thus $C$ is smooth at $p$ for any $p\notin \Sing(\Ff)$. Therefore, $\omega:=df$ generates $C$ around $p$, and we can pick $g=h:=1$ and $\eta:=0$ as in Definition \ref{defn: z and cs of foliation}, which implies (3).

(4) follows from \cite[Chapter 3, p39]{Bru00} if $p$ is a reduced singularity of $\Ff$. By Lemma \ref{lem: type of pd sing}, we may assume that $p$ is a Poincar\'{e}-Dulac singularity of \textbf{\rm \textbf{Type B}}, and suppose that $\Ff$ is generated by the $1$-form 
$$(nw+z^n)dz-zdw$$ 
around $p$ for the local coordinates $z$ and $w$. Since $C$ is $\Ff$-invariant, by \cite[Fact I.2.4]{McQ08}, $C=\{z=0\}$. In particular, we may pick $g:=1$ and $h:=nw+z^n$ as in Definition \ref{defn: z and cs of foliation}, which implies that $Z(\Ff,C,p)=1$.
\end{proof}

\begin{defn}
Let $X\in\mathcal{S}_{cyc}$ be a surface and $C$ a curve on $X$ such that each irreducible component of $C$ is $\Ff$-invariant. By Lemma \ref{lem: basic properties of z and cs}, we may define
$$Z(\Ff,C):=\sum_{p\in C}Z(\Ff,C,p)$$
and
$$CS(\Ff,C):=\sum_{p\in C}CS(\Ff,C,p).$$ 
\end{defn}

\subsection{$\Ff$-chains} 
\begin{defn}\label{defn: string on surface} Let $X$ be a normal projective surface and $C$ a curve on $X$. $C$ is called a \emph{string}, if
\begin{itemize}
\item $C=\cup_{i=1}^{n}C_i$ for some integer $n>0$, where each $C_i$ is a smooth rational curve, and
\item for any $1\leq i,j \leq n$, 
\begin{itemize}
    \item $C_i\cdot C_j=0$ if $|i-j|\geq 2$,
    \item $C_i\cdot C_j=1$ if $|i-j|=1$, and 
    \item $C_i^2<0$.  
\end{itemize}
\end{itemize} 
In this case, we define the \emph{intersection matrix} of $C$ to be the $n\times n$ matrix $\{(C_i\cdot C_j)\}_{n\times n}$, and is denoted by $||C||$.
\end{defn}

\begin{defn}\label{defn: F chain}  Let $X$ be a normal projective surface and $\Ff$ a foliation on $X$. An $\Ff$-\emph{chain} is a curve $C\subset X$, such that
\begin{itemize}
\item $C=\cup_{i=1}^{n}C_i$ is a string for some integer $n>0$,
\item each irreducible component $C_i$ of $C$ is $\Ff$-invariant,
\item $\Ff$ has at most reduced singularities along $C$, 
\item $Z(\Ff, C_1)=1$,
\item for any $i\geq 2$, $Z(\Ff, C_i)=2$, 
\item $C_1\cap\Sing(X)$ contains at most $1$ point. Moreover, if $C_1\cap\Sing(X)\not=\emptyset$, then $C_1\cap\Sing(X)$ is a cyclic quotient singularity, and
\item $C_i\cap \Sing(X)=\emptyset$ for any $i\geq 2$.
\end{itemize}
$C_1$ is called the \emph{initial curve} of $C$ and $C_n$ is called the \emph{last curve} of $C$. 
\end{defn}

\begin{rem}\label{rem: abuse notation fchain}
In the rest of the paper, if there's no confusion, when we say that ``$C=\cup_{i=1}^{n}C_i$ is an $\Ff$-chain", we always additionally assume that 
\begin{itemize}
\item $n>0$ is an integer,
    \item $C_1$ is the initial curve of $C$,
    \item $C_n$ is the last curve of $C$, and
    \item $C_i\cdot C_{i+1}=1$ for any $1\leq i\leq n-1$.
\end{itemize} 
\end{rem}

\begin{defn}[Maximal $\Ff$-chain]\label{defn: maximal F chain} Let $X$ be a normal projective surface, $\Ff$ a foliation on $X$, $n>0$ an integer and $V$ a set. An $\Ff$-chain $C=\cup_{i=1}^nC_i\in V$ is called \emph{maximal in $V$} if there does not exist a smooth rational curve $C_{n+1}\not\subset\Supp C$ such that $C'=\cup_{i=1}^{n+1}C_i$ is an $\Ff$-chain and $C'\in V$.
\end{defn}

\begin{lem}[Singularities along $\Ff$-chains]\label{lem: sing fchain basic}
Let $X$ be a normal projective surface, $\Ff$ a foliation on $X$, $n>0$ an integer and $C=\cup_{i=1}^nC_i$ an $\Ff$-chain on $X$. Let $p_i:=C_{i}\cap C_{i+1}$ for any $1\leq i\leq n-1$. Then there exists a closed point $p_n\in C_n$, such that
 \begin{enumerate}
 \item $p_n\not=p_{n-1}$, 
\item $C\cap\Sing(\Ff)=\{p_1,\dots, p_{n}\}.$ In particular, $p_i\not\in\Sing(X)$ for any $1\leq i\leq n$, 
\item for any $1\leq i\leq n$, $Z(\Ff,C_i,p_i)=1$, and
\item for any $1\leq i\leq n-1$, $Z(\Ff,C_{i+1},p_i)=1$.
 \end{enumerate}
\end{lem}
\begin{proof}
 The existence of $p_n$ satisfying (1)(2) follows from {\cite[Chapter 8, Definition 1]{Bru00}}. By Lemma \ref{lem: basic properties of z and cs}(4), such $p_n$ satisfies (3)(4).
\end{proof}

\begin{defn}[Tail of an $\Ff$-chain]\label{defn: tail of F chain} Let $X$ be a normal projective surface, $\Ff$ a foliation on $X$, $n>0$ an integer, and $C=\cup_{i=1}^{n}C_i$ an $\Ff$-chain. A \emph{tail} of $C$ is a curve $C_{n+1}\subset X$, such that 
\begin{itemize}
\item $C_{n+1}$ is an irreducible $\Ff$-invariant curve, 
\item $C_{n+1}\not\subset\Supp C$, and 
\item $C_{n+1}\cap C_n=p\in\Sing(\Ff)$. 
\end{itemize}
\end{defn}

\begin{lem}\label{lem: uniqueness of tail}
Let $X$ be a normal projective surface, $\Ff$ a foliation on $X$, and $C=\cup_{i=1}^{n}C_i$ an $\Ff$-chain. Then there exists at most $1$ tail of $C$.
\end{lem}

\begin{proof}
Suppose that $C_{n+1}$ and $C_{n+1}'$ are two tails of $C$, such that $C_{n+1}\cap C_{n}=q_1$ and $C_{n+1}'\cap C_{n}=q_2$ where $q_1$ and $q_2$ are closed points on $C_n$. If $q_1=q_2$, then there are at least three different $\Ff$-invariant curves passing through $q_1$, which is not possible since $q_1\in C_n$ is a reduced singularity. If $q_1\not=q_2$, then $q_1,q_2\in\Sing(\Ff)\cap C_n$, which contradicts to Lemma \ref{lem: sing fchain basic}(2).
\end{proof}

\begin{rem}\label{rem: the tail of c}
Let $X$ be a normal projective surface, $\Ff$ a foliation on $X$, and $C$ an $\Ff$-chain. By Lemma \ref{lem: uniqueness of tail}, we may call a tail of $C$ as \textbf{the} tail of $C$.
\end{rem}

At the end of this subsection, we introduce the notation of $(K_{\Ff}+\Delta)$-chains, which plays an important role in the construction of the minimal model program of foliated surface pairs.

\begin{defn}[$(K_{\Ff}+\Delta)$-chain and artificial chain]\label{defn: artificial chain} Let $(X,\Ff,\Delta)$ be a foliated triple such that $\dim X=2$, $n>0$ an integer, and $C=\cup_{i=1}^{n}C_i$ an $\Ff$-chain on $X$. We say that $C$ is a $(K_{\Ff}+\Delta)$-chain if there is a sequence of birational morphisms between foliated triples
$$(X_0:=X,\Ff_0:=\Ff,\Delta_0:=\Delta)\xrightarrow{\phi_1} (X_1,\Ff_1,\Delta_1)\xrightarrow{\phi_2}\dots\xrightarrow{\phi_n}(X_n,\Ff_n,\Delta_n)$$
satisfying the following. For every $1\leq i\leq n$, let $C_i'$ be the strict transform of $C_i$ on $X_{i-1}$, then
\begin{itemize}
    \item $\phi_i$ is the contraction of $C_i'$, 
    \item $(K_{\Ff_{i-1}}+\Delta_{i-1})\cdot C_i'<0$, and
    \item $(C_i')^2<0$.
\end{itemize}
Moreover,
\begin{itemize}
    \item assume that $p=(\phi_n\circ\dots\circ\phi_1)_*C$. If $p$ is a smooth point of $X_n$, then $C$ is called a \emph{$(K_\Ff+\Delta)$-artificial chain}.
    \item $C$ is called a \emph{maximal $(K_\Ff+\Delta)$-chain} if it is a maximal $\Ff$-chain in the set of all the $(K_\Ff+\Delta)$-chains.
    \item $C$ is called a \emph{maximal $(K_\Ff+\Delta)$-artificial chain} if it is a maximal $\Ff$-chain in the set of all the $(K_\Ff+\Delta)$-artificial chains. 
\end{itemize}
\end{defn}

\subsection{Family of foliations} 

\begin{defn}[Smooth family of foliations]\label{defn: smooth family of foliations} Let $X$ be a smooth threefold, $\Ff$ a rank 1 foliation on $X$, and $\pi: X\rightarrow T$ a smooth morphism to a curve. For any closed point $t\in T$, let $X_t:=\pi^{-1}(t)$ be the fiber of $\pi$ over $t$. $(X_t, \Ff_t)|_{t\in T}$ is called \emph{a smooth family of (canonical) foliations of surfaces} associated with $\pi$, if 
\begin{itemize}
\item $\Ff$ is tangent to the fibers of $\pi$, and
\item $\Sing(\Ff)$ is a pure codimension $2$ set of $X$, 
\end{itemize}
and for any closed point $t\in T$,
\begin{itemize}
\item  $X_t$ is a smooth projective surface,  
\item $\Sing(\Ff)\cap X_t$ is a finite set, and
\item $\Ff_t=\Ff|_{X_t}$ is a rank 1 (canonical) foliation on $X_t$.
\end{itemize}
\end{defn}

\subsection{Miscellaneous birational geometry definitions}

\begin{defn}[Algebraic Euler characteristic]\label{defn: euler char of curve} Let $X$ be a normal projective surface. For any curve $C\subset X$, we define the algebraic Euler characteristic $\chi(C)$ to be 
$$\chi(C):=-K_X\cdot C-C^2.$$ \end{defn}

\begin{defn}[Stable base locus]\label{defn: stable base locus} Let $X$ be a normal projective variety and $D$ a $\mathbb Q$-Cartier $\mathbb Q$-divisor on $X$. The \emph{stable base locus} of $D$ is
$$\cap_{m>0}\Bs(|mD|),$$
and is denoted by $\textbf{\rm\textbf{B}}(D)$. Let $A$ be an ample divisor on $X$ and $0<\epsilon\ll 1$ a rational number. We define 
$\textbf{\rm\textbf{B}}_{+}(D):=\textbf{\rm\textbf{B}}(D-\epsilon A)$
and
$\textbf{\rm\textbf{B}}_{-}(D):=\textbf{\rm\textbf{B}}(D+\epsilon A).$

 $\textbf{\rm\textbf{B}}_{+}(D)$ and $\textbf{\rm\textbf{B}}_{-}(D)$ are well-defined (cf. ~\cite[Chapter 10.3]{Laz04},~\cite{Nak00}.)
\end{defn}

\begin{defn}[Volume]\label{defn: volume}
Let $X$ be a normal projective variety and $D$ an $\mathbb R$-Cartier $\mathbb R$-divisor on $X$. We define the volume of $D$ to be
$$\vol(D):=\lim_{m\rightarrow+\infty}\sup\frac{n!h^0(X,mD)}{m^n}.$$
\end{defn}

\begin{defn}[Null]\label{defn: null curves} Let $X$ be a normal projective variety and $P$ a nef $\mathbb R$-divisor on $X$. We define 
$$\Null(P):=\{C|C\text{ is an irreducible curve on }X, P\cdot C=0\}.$$ \end{defn}

\begin{defn}[Restriction of linear system]\label{defn: restriction linear system} Let $X$ be a smooth projective variety, $\pi: X\rightarrow T$ a smooth morphism, and $|D|$ a linear system on $X$. We define
$$|D|_t:=|D|_{\pi^{-1}(t)}.$$
\end{defn}

\section{The classical theory}

In this section, we review several well-known results, and give generalizations for some of them. Most of the lemmas are classical statements of foliated surfaces: see \cite{Bru00}, \cite{Bru02}, and \cite{McQ08} for more details.

\subsection{Basic formulas of invariants of foliated surfaces} 
\begin{lem}[Intersection numbers of non-$\Ff$-invariant curves]\label{lem: tang formula} Let $X\in\mathcal{S}_{cyc}$ be a surface, $\Ff$ a foliation on $X$ and $C$ a curve on $X$ such that each irreducible component of $C$ is \textbf{not} $\Ff$-invariant. Then
\begin{enumerate}
    \item $K_\Ff\cdot C=-C^2+\tang(\Ff,C)$, and
    \item $(K_\Ff+C)\cdot C\geq 0$.
\end{enumerate}
\end{lem}

\begin{proof} (1) follows from ~\cite[Chapter 2, Proposition 2]{Bru00} and ~\cite[Chapter 2, p5]{Bru02}. Since $\tang(\Ff,C)\geq 0$, (1) implies (2).\end{proof}

\begin{lem}[Intersection numbers of $\Ff$-invariant curves, {~\cite[Chapter 2, Proposition 3]{Bru00}, ~\cite[Chapter 3, Theorem 2]{Bru00},~\cite[Chapter 2, p5]{Bru02}}]\label{lem: z and cs formula} Let $X\in\mathcal{S}_{cyc}$ be a surface, $\Ff$ a foliation on $X$ and $C$ a curve on $X$ such that each irreducible component of $C$ is $\Ff$-invariant. Then
\begin{enumerate}
\item $K_\Ff\cdot C=-\chi(C)+Z(\Ff,C)$, and
\item $C^2= CS(\Ff,C)$.
\end{enumerate}
\end{lem}

\begin{lem}[Eigenvalues and invariants of foliations]\label{lem: eigenvalue and cs and z formula} Assume that 
\begin{itemize}
    \item $X\in\mathcal{S}_{cyc}$ is a surface, 
    \item $\Ff$ is a canonical foliation on $X$, 
    \item $p\in\Sing(\Ff)$ is a closed point,
    \item  $C$ is an irreducible $\Ff$-invariant curve such that $p\in C$, and
    \item $\lambda:=\lambda(\Ff,p)$.
\end{itemize}
Then we have the following.
\begin{enumerate}
\item If $\lambda=0$, then
\begin{enumerate}
\item either $CS(\Ff,C,p)=0$ and $Z(\Ff,C,p)=1$, or 
\item $Z(\Ff,C,p)\geq 2$.
\end{enumerate}
\item If $\lambda\not=0$, then
\begin{enumerate}
\item either $CS(\Ff,C,p)=\lambda$, or 
\item $CS(\Ff,C,p)=\frac{1}{\lambda}$.
\end{enumerate}
\item If $p$ is reduced, then 
\begin{enumerate}
\item $\lambda\not=0$, 
\item there are exactly two irreducible $\Ff$-invariant curves $C$ and $C'$ passing through $p$, and
\item $\{CS(\Ff,C,p),CS(\Ff,C',p)\}=\{\lambda,\frac{1}{\lambda}\}.$
\end{enumerate}
\item If $p$ is Poincar\'{e}-Dulac, then 
\begin{enumerate}
\item $C$ is the unique irreducible $\Ff$-invariant curve passing through $p$, and 
\item $CS(\Ff,C,p)>0$.
\end{enumerate}
\end{enumerate}
\end{lem}

\begin{proof} Since $\Ff$ is canonical, by Lemma \ref{lem: either type b pd or reduced}, $p$ is either a reduced singularity or a Poincar\'{e}-Dulac singularity of  \textbf{\rm \textbf{Type B}}.  If $p$ is reduced, the proof follows from \cite[Chapter 3, Theorem 2]{Bru00}, so we may assume that $p$ is a Poincar\'{e}-Dulac singularity of  \textbf{\rm \textbf{Type B}}. In particular, we may assume that $\lambda\not=0$. Thus we only need to prove (2) and (4).

(4.a) follows from \cite[Fact I.2.4]{McQ08}. Possibly replacing $\lambda$ with $\frac{1}{\lambda}$, we may assume that $\lambda\in\mathbb N^+$. By Lemma \ref{lem: type of pd sing}, we may assume that $\Ff$ is generated by the $1$-form
$$(\lambda w+z^\lambda)dz-zdw$$
around $p$. By \cite[Fact I.2.4]{McQ08}, $C=\{z=0\}$. Thus
$$CS(\Ff, C, p)=\Res_0\{-\frac{1}{-\lambda w}dw|_{z=0}\}=\frac{1}{\lambda},$$
which implies (2) and (4.b).\end{proof}

We end this subsection by stating two lemmas that characterize the change of invariants for smooth foliated surfaces under blow-downs and blow-ups.
\begin{lem}\label{lem: contraction changes z}  Let $X$ be a smooth projective surface, $\Ff$ a canonical foliation on $X$, and $C$ a smooth rational $\Ff$-invariant curve on $X$ such that 
\begin{itemize}
\item $C^2=-1$, 
\item $Z(\Ff,C)=1$ or $2$, and 
\item $\Ff$ only has reduced singularities along $C$.
\end{itemize}
Suppose that $\nu: X\rightarrow X'$ is the contraction of $C$. Then 
\begin{enumerate}
\item $\nu$ is a blow-down to a smooth surface, 
\item $\nu_*\Ff$ is a canonical foliation,
\item if $C$ contains exactly one point of $\Sing(\Ff)$ and $Z(\Ff,C)=1$, then
\begin{itemize}
\item $\nu$ is the blow-down to a smooth point of $\nu_*\Ff$, and
\item for any $\Ff$-invariant curve $C'$ on $X$ such that $C'\cap C$ is a reduced singularity of $\Ff$, $Z(\nu_*\Ff,\nu_*C')=Z(\Ff,C')-1$, 
\end{itemize}
and
\item if $C$ contains exactly two points of $\Sing(\Ff)$ and $Z(\Ff,C)=2$, 
\begin{itemize}
\item $\nu$ is the blow-down to a reduced singularity of $\nu_*\Ff$, and
\item for any $\Ff$-invariant curve $C'$ on $X$ such that $C'\cap C$ is a reduced singularity of $\Ff$, $Z(\nu_*\Ff,\nu_*C')=Z(\Ff,C')$.
\end{itemize}
\end{enumerate}
\end{lem}
\begin{proof} (1) is Castelnuovo's theorem. Since $\Ff$ only has reduced singularities along $C$, \cite[Chapter 5, p72]{Bru00} shows that $\nu_*\Ff$ is canonical, which implies (2). (3)(4) follow from \cite[Chapter 5, p72]{Bru00}. \end{proof}

\begin{lem}[{\cite[Chapter 8, p114-115]{Bru00}}]\label{lem: blow up changes z}
Assume that
\begin{itemize}
    \item $X$ be a smooth projective surface,
    \item $\Ff$ is a canonical foliation on $X$,
    \item $p\in X$ is a closed point,
    \item $C$ is an $\Ff$-invariant curve such that $p\in C$,
    \item  $\pi: \bar{X}\rightarrow X$ is the blow-up of $p$, and
    \item $E$ is the exceptional curve of $\pi$.
\end{itemize}
Then
\begin{enumerate}
\item $E^2=-1$, 
\item if $p\notin \Sing(\Ff)$, then 
\begin{enumerate}
\item $K_{\pi^*{\Ff}}=\pi^*{K_{\Ff}}+E$, 
\item $Z(\Ff,E)=1$, 
\item $E$ contains exactly one point of $\Sing(\pi^*\Ff)$ which is a reduced singularity, and
\item $Z(\Ff,C)=Z(\pi^*\Ff,\pi^{-1}_*C)-1,$
\end{enumerate}
and
\item if $p\in \Sing(\Ff)$ and $p$ is reduced, then
\begin{enumerate}
\item $K_{\pi^*{\Ff}}=\pi^*{K_{\Ff}}$, 
\item $Z(\Ff,E)=2$, a
\item $E$ contains exactly two points of $\Sing(\pi^*\Ff)$ which are both reduced, and
\item $Z(\Ff,C)=Z(\pi^*\Ff,\pi^{-1}_*C).$
\end{enumerate}
\end{enumerate}
\end{lem}

\subsection{Contraction of $K_{\Ff}$-chains}

\begin{lem}[Contraction of a $K_{\Ff}$-negative $\Ff$-invariant curve]\label{lem: contraction of negative extremal ray foliation} Let $X\in\mathcal{S}_{cyc}$ be a surface, $\Ff$ a canonical foliation on $X$, and $C$ an irreducible $\Ff$-invariant curve such that $K_\Ff\cdot C<0$ and $C^2<0$. Then there exists a contraction $\nu: X\rightarrow X'$ of $C$, such that
\begin{enumerate}
\item $\nu_*\Ff$ is a canonical foliation,
\item $\nu_*C$ is either a smooth point of $X'$, or a cyclic quotient singularity of $X'$, and
\item $\nu_*C\notin \Sing(\nu_*\Ff)$.
\end{enumerate}
\end{lem}

\begin{proof}
By Lemma \ref{lem: z and cs formula} and following the steps as in \cite[Chapter 4, p10]{Bru02}, $C$ is a $K_\Ff$-negative extremal ray and it is contractible. The lemma follows from \cite[Chapter 4, p10]{Bru02}.
\end{proof}

\begin{lem}[Contraction of a $K_{\Ff}$-chain]\label{lem: F chain and KF chain}
Assume that 
\begin{itemize}
\item $X\in\mathcal{S}_{cyc}$ is a surface, 
\item $\Ff$ is a canonical foliation on $X$,
\item $n>0$ is an integer,
\item $C=\cup_{i=1}^{n}C_i$ is a $K_\Ff$-chain,
\item $C_{n+1}$ is the tail of $C$ if it exists, and
\item there exists a sequence of contractions 
$$(X_0:=X,\Ff_0:=\Ff,0)\xrightarrow{\phi_1} (X_1,\Ff_1,0)\xrightarrow{\phi_2}\dots\xrightarrow{\phi_n}(X_n,\Ff_n,0),$$
such that for any $1\leq i\leq n$, $\phi_i$ is the contraction of the strict transform of $C_i$ on $X_{i-1}$.
\end{itemize}
Then
\begin{enumerate}
    \item for every $1\leq i\leq n$, $X_i\in\mathcal{S}_{cyc}$,
    \item  for every $1\leq i\leq n$, $\Ff_i$ is a canonical foliation,
    \item  for every $1\leq i\leq n$, $(\phi_i\circ\phi_{i-1}\circ\dots\phi_1)_*(\cup_{j=1}^iC_j)\not\in\Sing(\Ff_i)$,
    \item for every $1\leq i\leq n-1$, $(\phi_i\circ\phi_{i-1}\circ\dots\phi_1)_*C$ is a $K_{\Ff_i}$-chain, and
    \item for every $1\leq i\leq n$, $Z(\Ff_i,(\phi_i\circ\phi_{i-1}\circ\dots\phi_1)_*C_{i+1})=Z(\Ff,C_{i+1})-1$.
\end{enumerate}
\end{lem}

\begin{proof}  Since $C$ is a $\Ff$-chain, $K_\Ff\cdot C_1<0$ and $C_1^2<0$, and there are closed points $p_1,\dots,p_n$ on $X$ such that
\begin{itemize}
    \item  $p_1,\dots, p_n$ are the only singularities of $\Ff$ along $C$,
    \item $p_i=C_i\cap C_{i+1}$ for each $i$, and
    \item $Z(\Ff,C_i,p_i)=Z(\Ff,C_{i+1},p_i)=1$ for each $i$.
\end{itemize}
By Lemma \ref{lem: contraction of negative extremal ray foliation}(3) and Lemma \ref{lem: basic properties of z and cs}(3), 
$$Z(\Ff_1,(\phi_1)_*C_2)=1$$ and 
$$Z(\Ff_1,(\phi_1)_*C_k)=2$$
for any $3\leq k\leq n$. Thus $(\phi_1)_*C$ is a $\Ff_1$-chain. Since $C$ is a $K_{\Ff}$-chain, $(\phi_1)_*C$ is a $K_{\Ff_1}$-chain with length $n-1$. The lemma follows from the induction on the length of $C$.
\end{proof}

\subsection{A detailed description of $\Ff$-chains} To generalize the theory of the minimal model program for foliated surfaces to the minimal model program for foliated surface pairs, in this subsection we give a more detailed description on $K_{\Ff}$-chains and $(K_{\Ff}+\Delta)$-chains.

\begin{lem}\label{lem: Intersection surface} Let $X$ be a normal projective surface, $E$ an irreducible curve on $X$ and $f: X\rightarrow X'$ the contraction of $E$. Then for any $\mathbb R$-divisor $D$ on $X$ and any curve $C$ on $X$, 
$$f_*D\cdot f_*C=D\cdot C+\frac{(C\cdot E)}{(-E^2)}(D\cdot E).$$
\end{lem}

\begin{proof}
Let $D':=f_*D$ and $C':=f_*C$, and suppose that $f^*D'=D+aE$ for some real number $a$. Since
$f^*D'\cdot E=0,$ we have $a=\frac{D\cdot E}{-E^2}.$ By the projection formula,
$$f_*D\cdot f_*C=D'\cdot C'=f^*D'\cdot C=(D+aE)\cdot C=D\cdot C+\frac{(C\cdot E)}{(-E^2)}(D\cdot E).$$
\end{proof}

The next lemma shows that maximal $(K_{\Ff}+\Delta)$-chains do not intersect.
\begin{lem}\label{lem: different F chains maximal} Let $(X,\Ff,\Delta)$ be a foliated triple, such that 
\begin{itemize}
    \item $X\in\mathcal{S}_{cyc}$,
    \item $\Ff$ is a canonical, and
    \item $K_\Ff+\Delta$ is pseudo-effective.
\end{itemize} 
Assume that 
\begin{itemize}
    \item $V$ is a set of $(K_{\Ff}+\Delta)$-chains,
    \item $n,m>0$ are two integers,
    \item $C=\cup_{i=1}^nC_i\in V$ and $C'=\cup_{i=1}^mC_i'\in V$ are two different $(K_{\Ff}+\Delta)$-chains such that $C\cap C'\not=\emptyset$.
\end{itemize} 
then 
\begin{enumerate}
\item either $C\subset C'$ or $C'\subset C$, and in particular,
\item either $C$ is not maximal, or $C'$ is not maximal.
\end{enumerate}
\end{lem}

\begin{proof} 
By our assumptions, there are closed points $p_1,\dots,p_n$ and $q_1,\dots, q_m$ on $X$, such that
\begin{itemize}
\item $p_1,\dots, p_n$ are the only singularities of $\Ff$ along $C$,
\item $q_1,\dots, q_m$ are the only singularities of $\Ff$ along $C'$,
\item for any $1\leq i\leq n-1$, $p_i=C_i\cap C_{i+1}$, and
\item for any $1\leq i\leq m-1$, $q_i=C'_i\cap C'_{i+1}$.
\end{itemize}
Without loss of generality, we may assume that $n\geq m$. There are two cases:\vspace{2mm}
\noindent\textbf{Case A}. $m\geq 2$, and $\{q_1,\dots, q_{m-1}\}\subset \{p_1,\dots, p_{n-1}\}$.

In this case, by Lemma \ref{lem: eigenvalue and cs and z formula}(3), $C'$ is the union of all the irreducible $\Ff$-invariant curves passing through $q_1,\dots, q_{m-1}$ and $C$ is the union of all the irreducible $\Ff$-invariant curves passing through $p_1,\dots p_{n-1}$. In particular, $C'\subset C$, and the Lemma follows.\vspace{2mm}

\noindent\textbf{Case B}. 
\begin{itemize}
\item Either $m=1$, or 
\item $m\geq 2$, and there exists an integer $1\leq i\leq m-1$ such that $q_i\notin\{p_1,\dots,p_{n-1}\}$.
\end{itemize}
\begin{claim}\label{claim: cr intersects c maximal chain}
In \textbf{\rm \textbf{Case B}} as above, we may assume that there exists an integer $1\leq r\leq m$, such that 
\begin{itemize}
\item $C_r'\cap C_n\not=\emptyset$, 
\item $C_r'$ is not an irreducible component of $C$, and 
\item $C_1'\dots,C_{r-1}'$ do not intersect $C$. 
\end{itemize}
\end{claim}
\begin{proof}[Proof of Claim \ref{claim: cr intersects c maximal chain}] There are two cases: $m=1$, or $m\geq 2$.\vspace{2mm}

\noindent\textbf{Case 1}. $m=1$. In this case, $C'=C_1'$. If $C_1'\subset C$ then we are done. Otherwise, by Lemma \ref{lem: eigenvalue and cs and z formula}(3), $q_1=p_n$. Thus $C_1'$ is not an irreducible component of $C$, and $C_1'\cap C_n\not=\emptyset$. We may take $r:=1$.\vspace{2mm}

\noindent\textbf{Case 2}. $m\geq 2$. In this case, either $q_1\in C$, or $q_1\not\in C$.\vspace{2mm}

\noindent\textbf{Case 2.1} $q_1\in C$. By Lemma \ref{lem: eigenvalue and cs and z formula}(3), there are only two possible cases:\vspace{2mm}

\noindent\textbf{Case 2.1.1} $q_1=p_1$. By inductively applying Lemma \ref{lem: eigenvalue and cs and z formula}(3), we deduce that $p_2=q_2,\dots, p_{m}=q_{m}$, which contradicts to our assumptions.\vspace{2mm}

\noindent\textbf{Case 2.1.2} $q_1=p_n$ and $C_1'\cap C=p_n$. In this case we may let $r:=1$.\vspace{2mm}

\noindent\textbf{Case 2.2} $q_1\not\in C$. Thus $C_1'$ does not intersect $C$. By Lemma \ref{lem: eigenvalue and cs and z formula}(3), we may let $r\geq 2$ be the integer such that $q_r\in C$ but $q_1,\dots, q_{r-1}\notin C$.
\end{proof}
\noindent\textit{Proof of Lemma \ref{lem: different F chains maximal} continued}. Pick $r>0$ as in Claim \ref{claim: cr intersects c maximal chain}. By Lemma \ref{lem: eigenvalue and cs and z formula}(3), $C_r'\cap C=p_n$.

Since $C'$ is a $(K_\Ff+\Delta)$-chain, $\cup_{i=1}^rC_i'$ is a $(K_\Ff+\Delta)$-chain. Thus there exists a birational contraction $\rho: X\rightarrow X'$ of $\cup_{i=1}^rC_i'$ satisfying the following.
\begin{itemize}
\item For any $1\leq i\leq n-1$, since  $C_i$ does not intersects $C'$, we have
\begin{equation}\label{equ: contract first maximal chain ci}
(K_{\rho_*\Ff}+\rho_*\Delta)\cdot \rho_*C_{i}=(K_{\Ff}+\Delta)\cdot C_{i}<0.\tag{3.9.1}
\end{equation}
\item By Lemma \ref{lem: Intersection surface},
\begin{equation}\label{equ: contract first maximal chain cn}
(K_{\rho_*\Ff}+\rho_*\Delta)\cdot \rho_*C_{n}<(K_{\Ff}+\Delta)\cdot C_{n}<0. \tag{3.9.2}
\end{equation}
\item For any $1\leq i\leq n$, since $K_\Ff+\Delta$ is pseudo-effective, by (3.9.1) and (3.9.2),
\begin{equation}\label{equ: contract first maximal chain ci2}
(\rho_*C_i)^2<0.\tag{3.9.3}
\end{equation}
\end{itemize}

If $n=1$, then $p_1$ is the only singularity of $\Ff$ along $C$. By Lemma \ref{lem: contraction of negative extremal ray foliation}(3), $\rho_*\Ff$ is smooth along $\rho_*C$. By Lemma \ref{lem: basic properties of z and cs}(3) and Lemma \ref{lem: z and cs formula}(2), $(\rho_*C)^2=0$. This contradicts to (3.9.3).

If $n\geq 2$, by (3.9.1), (3.9.2) and (3.9.3), $\cup_{i=1}^{n-1}\rho_*C_i$ is a $(K_{\rho_*\Ff}+\rho_*\Delta)$-chain. Let $\tau$ be the contraction of $\cup_{i=1}^{n-1}\rho_*C_i$. By Lemma \ref{lem: Intersection surface}, 
\begin{equation}\label{equation: section contraction maximal chain}
(K_{(\tau\circ\rho)_*\Ff}+(\tau\circ\rho)_*\Delta)\cdot (\tau\circ\rho)_*C_{n}\leq (K_{\rho_*\Ff}+\rho_*\Delta)\cdot \rho_*C_{n}<0.\tag{3.9.4}
\end{equation}
Since $K_\Ff+\Delta$ is pseudo-effective, we have
\begin{equation}\label{equation: section contraction maximal chain self intersection}
((\tau\circ\rho)_*C_{n})^2<0.\tag{3.9.5}
\end{equation}
By Lemma \ref{lem: F chain and KF chain}(3),  $(\tau\circ\rho)_*\Ff$ is smooth along $(\tau\circ\rho)_*C_n$. By Lemma \ref{lem: basic properties of z and cs}(3) and Lemma \ref{lem: z and cs formula}(2),  $((\tau\circ\rho)_*C_n)^2=0$. This contradicts to (3.9.5). \end{proof}

\begin{lem}\label{lem: contract -1curve in a chain}
Assume that 
\begin{itemize}
    \item $(X,\Ff,\Delta)$ is a foliated triple such that $X\in\mathcal{S}_{cyc}$ and $\Ff$ is canonical,
    \item $n>0$ is an integer, and
    \item $C=\cup_{k=1}^nC_k$ a $(K_{\Ff}+\Delta)$-chain.
\end{itemize} 
Then for any integer $1\leq i\leq n$ such that $C_i^2=-1$, there is a contraction $\phi: (X,\Ff,\Delta)\rightarrow (X',\Ff',\Delta')$, such that
\begin{enumerate}
    \item $\phi$ is the contraction of $C_i$,
    \item $X'\in\mathcal{S}_{cyc}$, and if $X$ is smooth, then $X'$ is smooth,
    \item $\Ff'$ is canonical,
    \item $\phi_*C$ is either a point or a $(K_{\Ff'}+\Delta')$-chain.
\end{enumerate}
\end{lem}
\begin{proof}
The existence of $\phi$ and (1)(2)(3) follow from Lemma \ref{lem: contraction changes z}, and we only left to prove (4). We may assume that $i\geq 2$, otherwise $C_i=C_1$ itself is a $(K_{\Ff}+\Delta)$-artificial chain and there is nothing left to prove. Since $\nu$ is an isomorphism near $\cup_{k=1}^{i-2}C_k$, possibly contracting $\cup_{k=1}^{i-2}C_k$ and replacing $X$, we may assume that $i=2$. Let $C_1':=\phi_*C_1$, $\nu: (X,\Ff,\Delta)\rightarrow (X'',\Ff'',\Delta'')$ be the contraction of $C_2$ and $C_2'':=\nu_*C_2$, then by Lemma \ref{lem: Intersection surface}, we have 
$$(K_{\Ff'}+\Delta')\cdot C_1'=(K_{\Ff}+\Delta)\cdot C_1+(K_{\Ff}+\Delta)\cdot C_2,$$
$$0>(K_{\Ff''}+\Delta'')\cdot C_2''=(K_{\Ff}+\Delta)\cdot C_2+\frac{1}{-C_1^2}(K_{\Ff}+\Delta)\cdot C_1,$$
$$(C_1')^2=C_1^2+1,$$
and
$$0>(C_2'')^2=-1+\frac{1}{-C_1^2}.$$
We deduce that $C_1^2<-1$, $(C_1')^2<0$ and $(K_{\Ff'}+\Delta')\cdot C_1'<0$, which implies (4).
\end{proof}

\begin{lem}\label{lem: contract artificial chain step by step}  Assume that
\begin{itemize}
    \item $(X,\Ff,\Delta)$ is a foliated triple such that $X$ is a smooth surface and $\Ff$ is canonical,
    \item  $n>0$ is  an integer, 
    \item $C=\cup_{k=1}^{n}C_k$ is a $(K_\Ff+\Delta)$-artificial chain, and
    \item  $\nu: X\rightarrow X'$ is the contraction of $C$.
\end{itemize}  Then there is a sequence of contractions
$$(X_0,\Ff_0,\Delta_0):=(X,\Ff,\Delta)\xrightarrow{\nu_1} (X_1,\Ff_1,\Delta_1)\xrightarrow{\nu_2}\dots\xrightarrow{\nu_n}(X_n=X',\Ff_n,\Delta_n),$$
such that for every $1\leq i\leq n$,
\begin{enumerate}
    \item $\nu_i$ is a blow-down of the strict transform of an irreducible component of $C$,
    \item $X_i$ is smooth,
    \item $\Ff_i$ is a canonical foliation, and
    \item $(\nu_i\circ\dots\circ\nu_1)_*C$ is a $(K_{\Ff_i}+\Delta_i)$-chain if $i\leq n-1$.
\end{enumerate}
\end{lem}

\begin{proof} Since $X$ is smooth and $C$ is artificial, by \cite[Chapter 8, p113]{Bru00}, there exists an integer $1\leq k\leq n$ such that $C_{k}^2=-1$. Let $\nu_1: (X,\Ff,\Delta)\rightarrow (X_1,\Ff_1,\Delta_1)$ be the contraction of $C_{k}$. By Lemma \ref{lem: contraction changes z}(1), $\nu_1$ is a blow-down, $X_1$ is smooth, and $\Ff_1$ is canonical. By Lemma \ref{lem: contract -1curve in a chain}, $(\nu_1)_*C$ is a $(K_{\Ff_1}+\Delta_1)$-chain. Thus $(\nu_1)_*C$ is a $(K_{\Ff_1}+\Delta_1)$-artificial chain with length $n-1$. Possibly replacing $X,\Ff,C$ and $\nu$ with $X_1,\Ff_1,(\nu_1)_*C$ and the contraction of $(\nu_1)_*C$ respectively, we finish the proof by applying the induction on $n$.
 \end{proof}
 
 \begin{lem}\label{lem: no continue -1 curve}
 Assume that
 \begin{itemize}
    \item $(X,\Ff,\Delta)$ is a foliated triple such that $X$ is a smooth surface and $\Ff$ is canonical,
    \item $n>0$ is an integer, 
    \item $C=\cup_{k=1}^{n}C_k$ is a $(K_\Ff+\Delta)$-chain,
    \item $1\leq i,j\leq n$ are two integers such that $i\not=j$, and
    \item $C_i^2=C_j^2=-1$,
    \end{itemize}
then $C_i$ and $C_j$ do not intersect.
\end{lem}
\begin{proof}
By Lemma \ref{lem: contraction changes z}, there is a contraction $f: (X,\Ff,\Delta)\rightarrow (X',\Ff',\Delta')$ of $C_i$, such that $X'$ is smooth and $\Ff'$ is canonical.  By Lemma \ref{lem: contraction changes z}(3)(4) and Lemma \ref{lem: Intersection surface}, $f_*C$ is a $(K_{\Ff'}+\Delta')$-chain. 

If $C_i$ and $C_j$ intersects, then $C_i\cdot C_j=1$. By  Lemma \ref{lem: Intersection surface}, $(f_*C_j)^2=C_j^2+1=0$, which is not possible.
\end{proof}

\begin{lem}[{\cite[Remarks III.1.3]{McQ08}}]\label{lem: contraction intersection numbers} Assume that
\begin{itemize}
    \item $X$ is a smooth projective surface,
    \item $\Ff$ is a canonical foliation on $X$, 
    \item $n>0$ is an integer, 
    \item $C=\cup_{i=1}^{n}C_i$ is a $K_\Ff$-chain,
    \item $\nu: X\rightarrow X'$ is the contraction of $C$, and
    \item $C_{n+1}$ is the tail of $C$,
\end{itemize}
then $K_{\nu_*\Ff}\cdot \nu_*C_{n+1}=K_\Ff\cdot C_{n+1}-\det(-||C||)^{-1}$.
\end{lem}

\subsection{Classical deformation theory of foliated surfaces}

\begin{lem}\label{lem: F chain deforms} Assume that
\begin{itemize}
    \item $(X_t,\Ff_t)|_{t\in T}$ is a smooth family of canonical foliations of surfaces associated with $\pi: X\rightarrow T$, 
    \item $s\in T$ is a closed point, and
    \item $C_s=\cup_{i=1}^{n}C_{i,s}$ is an $\Ff_s$-chain such that $C_{i,s}^2\leq -2$ for any $1\leq i\leq n$,
\end{itemize} 
then there exists an open neighborhood $U$ of $s$ and hypersurfaces $C_1,\dots, C_n\subset \pi^{-1}(U)$ satisfying the following. For any $1\leq i\leq n$ and any closed point $t\in U$, 
\begin{enumerate}
\item $C_i$ is transversal to $X_t$, and
\item  $\cup_{i=1}^{n}(C_i\cap X_t)$ is an $\Ff_t$-chain.
\end{enumerate}
\end{lem}

\begin{proof} Since $\Ff_s$ only has reduced singularities along $C_s$, the lemma follows from \cite[Lemme 3]{Bru01}.\end{proof}

\begin{lem}\label{lem: -1 curve deforms}
Assume that
\begin{itemize}
\item $X$ is a smooth threefold,
\item $\pi: X\rightarrow T$ is a smooth projective morphism to a curve,
\item $X_t:=\pi^{-1}(t)$ for any $t\in T$,
\item $s\in T$ is a closed point, and
\item $C_s$ is a smooth rational curve on $X_s$ such that $C_s^2=-1$,
\end{itemize}
then there exists an open neighborhood $U$ of $s$ in $T$, a hypersurface $E\subset\pi^{-1}(U)$, a birational contraction $\nu: \pi^{-1}(U)\rightarrow X'$ and a smooth morphism $\pi':X'\rightarrow U$ satisfying the following.
\begin{enumerate}
\item  $E$ is transveral to the fibers of $\pi$, 
\item $E|_{X_s}=C_s$, 
\item for every closed point $t\in U$, $C_t:=E|_{X_t}$ is a smooth rational curve such that $C_t^2=-1$, 
\item $\pi|_{\pi^{-1}(U)}=\pi'\circ\nu|_{\pi^{-1}(U)}$, and 
\item for every closed point $t\in U$, $\nu_t:=\nu|_{X_t}$ is the contraction of $C_t$.
\end{enumerate}
\end{lem}

\begin{proof} (1)(2)(3) follow from \cite{Kod63}. If $\kappa(X)=-\infty$, (4)(5) follow from \cite[Theorem 1.11]{FM94}. If $\kappa(X)\geq 0$, (4)(5) follow from \cite[Theorem 1.12]{FM94} and \cite[Theorem 1.16]{FM94}.\end{proof}

\subsection{Several technical results}

\begin{lem}\label{lem: stable base locus under pushforward} Let $X$ be a normal projective variety, $f: X\rightarrow Y$ a birational contraction, and $\Delta$ and $G$ two $\mathbb Q$-Cartier $\mathbb Q$-divisors on $X$, such that 
\begin{itemize}
\item $\Delta\geq 0$,
\item $\kappa(G)\geq 0$, and 
\item $\Delta\wedge\textbf{\rm \textbf{B}}(G)=0$. 
\end{itemize}
Then $(f_*\Delta)\wedge \textbf{\rm \textbf{B}}(f_*G)=0.$
\end{lem}
\begin{proof}
Let $m>0$ be a sufficiently divisible integer, and write $\Delta=\sum_{i=1}^na_iD_i$ into its irreducible components. By our assumptions, for any $ 1\leq i\leq n$, there exists $G_i\in |mG|$ such that $D_i\not\subset\Supp G_i$. In particular, $f_*D_i\not\subset\Supp f_*G_i$. Since $f$ is a birational contraction, $f_*G_i\in |mf_*G|$. Therefore 
$$0=\sum a_if_*D_i\wedge\Supp(f_*G_i)\geq \sum a_if_*D_i\wedge \textbf{B}(f_*G)=(f_*\Delta)\wedge \textbf{B}(f_*G)\geq 0.$$
\end{proof}

\begin{lem}\label{lem: aug stable base locus pushforward}
Let $X$ be a normal projective surface, $f: X\rightarrow Y$ a sequence of divisorial contractions of curves with negative self-intersections, and $\Delta$ and $G$ two $\mathbb Q$-Cartier $\mathbb Q$-divisors on $X$, such that 
\begin{itemize}
\item $\Delta\geq 0$,
\item $G$ is pseudo-effective (resp. $G$ is big), and
\item $\Delta\wedge\textbf{\rm \textbf{B}}_{-}(G)=0$ (resp. $\Delta\wedge\textbf{\rm \textbf{B}}_{+}(G)=0$).
\end{itemize}
Then $(f_*\Delta)\wedge \textbf{\rm \textbf{B}}_{-}(f_*G)=0$ (resp. $f_*G$ is big and $(f_*\Delta)\wedge\textbf{\rm \textbf{B}}_{+}(f_*G)=0$).\end{lem}

\begin{proof} Let $A$ be an ample divisor on $X$. By Lemma \ref{lem: Intersection surface}, $f_*A$ is an ample divisor on $Y$. 

If $G$ is pseudo-effective and $\Delta\wedge\textbf{\rm \textbf{B}}_{-}(G)=0$, then there exists $0<\epsilon\ll 1$, such that 
$$\textbf{B}_{-}(G)=\textbf{B}(G+\epsilon A) \text{ and }\textbf{B}_{-}(f_*G)=\textbf{B}(f_*G+\epsilon f_*A).$$ 
Thus by Lemma \ref{lem: stable base locus under pushforward},
$$(f_*\Delta)\wedge\textbf{B}_{-}(f_*G)= (f_*\Delta)\wedge\textbf{B}(f_*G+\epsilon f_*A)=(f_*\Delta)\wedge\textbf{B}(f_*(G+\epsilon A))=0.$$
If $G$ is big and $\Delta\wedge\textbf{\rm \textbf{B}}_{+}(G)=0$, then there exists $0<\epsilon\ll 1$, such that 
$$\textbf{B}_{+}(G)=\textbf{B}(G-\epsilon A) \text{ and } \textbf{B}_{+}(f_*G)=\textbf{B}(f_*G-\epsilon f_*A).$$ 
Thus by Lemma \ref{lem: stable base locus under pushforward},
$$(f_*\Delta)\wedge\textbf{B}_{+}(f_*G)= (f_*\Delta)\wedge\textbf{B}(f_*G-\epsilon f_*A)=(f_*\Delta)\wedge\textbf{B}(f_*(G-\epsilon A))=0.$$
\end{proof}

\begin{lem}[Algebraic version of the negativity lemma]\label{lem: num negativity lemma}
Let $n>0$ be an integer.
\begin{enumerate}
\item Let $\bm{E}=\{e_{i,j}\}_{n\times n}$ be a negative definite matrix, such that $e_{i,j}\geq 0$ if $i\not=j$. If $\bm{a}=(a_1,\dots, a_n)$ is a vector such that all the coefficients of $\bm{aE}$ are greater or equal to zero, then $a_i\leq 0$ for any $1\leq i\leq n$.
\item Let $X,Y$ be normal varieties and $f: Y\rightarrow X$ a birational contraction. Let $D$ be an $f$-nef $f$-exceptional $\mathbb R$-divisor, then $-D$ is effective.
\end{enumerate}
\end{lem}
\begin{proof} (2) is the usual negativity lemma (cf. \cite[Lemma 3.39]{KM98}). To prove (1), write $\bm{a}=\bm{b}-\bm{c}$ where $\bm{b}=(b_1,\dots, b_n)$ and $\bm{c}=(c_1,\dots, c_n)$ are uniquely determined vectors such that for any $1\leq i\leq n$, $b_i\geq 0, c_i\geq 0$, and $b_ic_i=0$. Since all the coefficients of $\bm{aE}$ are greater or equal to zero, we have
$$\bm{aEb}^T=\bm{bEb}^T-\bm{cEb}^T\geq 0.$$
Since $E$ is negative definite, $\bm{bEb}^T\leq 0$. Since $e_{i,j}\geq 0$ for any $i\not=j$, $\bm{cEb}^T\geq 0$. Thus $\bm{bEb}^T=0$, an in particular, $\bm{b}=\bm{0}$. Thus $a_i\leq 0$ for any $1\leq i\leq n$.
\end{proof}

\begin{lem}\label{lem: Zariski decomposition and Null} Let $X$ be a normal projective surface and $K$ a pseudo-effective $\mathbb Q$-Cartier $\mathbb Q$-divisor on $X$, such that $K=P+N$ is the Zariski decomposition of $K$ where $P$ is the positive part and $N$ is the negative part. Then 
\begin{enumerate}
\item $\textbf{\rm \textbf{B}}_{-}(K)\subset\Supp N$, and
\item if $K$ is big, then $\textbf{\rm \textbf{B}}_{+}(K)=\Null(P)$.
\end{enumerate}
\end{lem}

\begin{proof} The statements follow from \cite[Example 2.19]{ELMNP09} and \cite[Example 3.4]{ELMNP06}.\end{proof}

\section{The minimal model program of foliated surface pairs} In this section we study the minimal model program for foliated surface pairs satisfying the additional assumption
$$\Delta\wedge\textbf{\rm \textbf{B}}_{-}(K_\Ff+\Delta)=0.$$
We give the following definitions, which will greatly simplify our statements below:
\begin{defn}
We define $\mathcal{D}$ to be the set of all the foliated triples $(X,\Ff,\Delta)$, such that
\begin{itemize}
   \item $X\in\mathcal{S}_{cyc}$ is a surface,
    \item $\Ff$ is a canonical foliation,
        \item $K_\Ff+\Delta$ is pseudo-effective, and
         \item $\Delta\wedge\textbf{\rm \textbf{B}}_{-}(K_\Ff+\Delta)=0$.
\end{itemize}
We define $\mathcal{D}_0$ to be the set of all the foliated triples $(X,\Ff,\Delta)\in\mathcal{D}$, such that
\begin{itemize}
    \item $X$ is smooth, 
    \item $K_\Ff+\Delta$ is big, 
    \item $\Delta\wedge\textbf{\rm \textbf{B}}_{+}(K_\Ff+\Delta)=0$, and
    \item $X$ does not contain any artificial $(K_\Ff+\Delta)$-chain.
\end{itemize}
\end{defn}

First, we study the structure of $(K_{\Ff}+\Delta)$-negative extremal rays.
\begin{lem}\label{lem: property of negative extremal ray foliated pairs} Let $(X,\Ff,\Delta)\in\mathcal{D}$ be a foliated triple and $C$ an irreducible curve on $X$ such that $(K_\Ff+\Delta)\cdot C<0$. Then
\begin{enumerate}
\item $C$ is $\Ff$-invariant,
\item $C$ is a smooth rational curve,
\item $Z(\Ff,C)=1$,
\item $p=C\cap \Sing(\Ff)$, such that
\begin{itemize}
    \item $p$ is a reduced singularity, and
    \item $\lambda:=\lambda(\Ff,p)$ is a negative rational number,
\end{itemize}
and
\item $C\cap\Sing(X)$ contains at most $1$ point. 
\end{enumerate}
\end{lem}

\begin{proof} Since $\Delta \wedge \textbf{B}_{-}(K_\Ff+\Delta)=0$ and $(K_\Ff+\Delta)\cdot C<0$, $\Delta\cdot C\geq 0$. Thus $K_\Ff\cdot C<0$. Since $K_\Ff+\Delta$ is pseudo-effective, $C^2<0$. Thus (1) follows from Lemma \ref{lem: tang formula}.

By Lemma \ref{lem: basic properties of z and cs}(3)(4) and Lemma \ref{lem: z and cs formula}(1), 
$$0>K_\Ff\cdot C= -\chi(C)+Z(\Ff, C)\geq -\chi(C),$$
hence $\chi(C)>0$, which implies (2). 

Since $\chi(C)\leq 2$,
$$Z(\Ff,C)=K_\Ff\cdot C+\chi(C)<\chi(C)\leq 2.$$ 
Since $C^2<0$, by Lemma \ref{lem: z and cs formula}(2), $CS(\Ff,C)<0$. By Lemma \ref{lem: basic properties of z and cs}(3), $C\cap \Sing(\Ff)\not=\emptyset$. Thus $C$ passes through at least $1$ closed point of $\Sing(\Ff)$. By (2) and Lemma \ref{lem: basic properties of z and cs}(2)(4), $Z(\Ff,C)\geq 1$. Thus $C\cap\Sing(\Ff)$ contains exactly $1$ point and $Z(\Ff,C)=1$, which implies (3).

Suppose that $p=C\cap \Sing(\Ff)$.  If $\lambda=0$,  by Lemma \ref{lem: eigenvalue and cs and z formula}(1), either $CS(\Ff,C,p)=0$ or $Z(\Ff,C,p)\geq 2$, both of which are not possible. Thus $\lambda\not=0$. By Lemma \ref{lem: z and cs formula}(2) and Lemma \ref{lem: eigenvalue and cs and z formula}(2), $$C^2=CS(\Ff,C,p)\in\{\lambda,\frac{1}{\lambda}\}.$$ 
Thus $\lambda$ is a negative rational number. By Lemma \ref{lem: eigenvalue and cs and z formula}(4) and Lemma \ref{lem: either type b pd or reduced}, $p$ is a reduced singularity, which implies (4).

Suppose that $C\cap\Sing(X)=\{p_1,\dots,p_k\}$ where each $p_i$ is a cyclic quotient singularity of order $r_i\geq 2$, then 
$$\chi(C)=2+\sum_{i=1}^k (\frac{1}{r_i}-1).$$ 
Since
$$0>-\chi(C)+Z(\Ff,C)=-\chi(C)+1,$$
$\chi(C)>1$, which implies that $k=1$, and we deduce (5).\end{proof}

\begin{lem}\label{lem: kfb chain is kf chain} Let $(X,\Ff,\Delta)\in\mathcal{D}$ be a foliated triple, then 
\begin{enumerate}
    \item every $(K_{\Ff}+\Delta)$-negative extremal ray is a $K_{\Ff}$-negative extremal ray,
    \item every $(K_{\Ff}+\Delta)$-chain is a $K_{\Ff}$-chain, and
    \item every partial $(K_{\Ff}+\Delta)$-MMP is a partial $K_{\Ff}$-MMP.
\end{enumerate}
\end{lem}

\begin{proof}
Let $C$ be an irreducible curve on $X$ such that $(K_\Ff+\Delta)\cdot C<0$. Since $\Delta\wedge\textbf{\rm \textbf{B}}_{-}(K_\Ff+\Delta)=0,$ $\Delta\cdot C\geq 0$. Thus $K_\Ff\cdot  C<0$, which implies that every $(K_{\Ff}+\Delta)$-negative extremal ray is a $K_{\Ff}$-negative extremal ray. By applying Lemma \ref{lem: aug stable base locus pushforward}, every $(K_\Ff+\Delta)$-chain is a $K_\Ff$-chain and every partial $(K_{\Ff}+\Delta)$-MMP is a partial $K_{\Ff}$-MMP.
\end{proof}

The next easy lemma is an important ingredient in the construction of the minimal model program for foliated surface pairs:
\begin{lem}\label{lem: contract reduce sing of f}
Let $(X,\Ff,\Delta)\in\mathcal{D}$ be a foliated triple and $f:(X,\Ff,\Delta)\rightarrow (X',\Ff',\Delta')$ a contraction, such that $f$ is 
\begin{itemize}
    \item either a contraction of a $(K_\Ff+\Delta)$-negative extremal ray, or
    \item a contraction of a $(K_\Ff+\Delta)$-chain.
\end{itemize}
Then
\begin{enumerate}
    \item $(X',\Ff',\Delta')\in\mathcal{D}$, and
    \item the number of closed points belonging to $\Sing(\Ff')$ is strictly less than the number of closed points belonging to $\Sing(\Ff)$.
\end{enumerate}
In particular, for every partial $(K_\Ff+\Delta)$-MMP $\rho: (X,\Ff,\Delta)\rightarrow (X'',\Ff'',\Delta'')$, $(X'',\Ff'',\Delta'')\in\mathcal{D}$.
\end{lem}
\begin{proof}
It follows from Lemma \ref{lem: contraction of negative extremal ray foliation}, Lemma \ref{lem: F chain and KF chain}(3), Lemma \ref{lem: aug stable base locus pushforward}, Lemma \ref{lem: property of negative extremal ray foliated pairs}(1)(2)(4), and Lemma \ref{lem: kfb chain is kf chain}.
\end{proof}

Now we show the existence of the minimal model program for foliated surface pairs:
\begin{lem}\label{lem: mmp foliated pairs weak version}Let $(X,\Ff,\Delta)\in\mathcal{D}$ be a foliated triple and $\rho: (X,\Ff,\Delta)\rightarrow (X',\Ff',\Delta')$ a partial $(K_\Ff+\Delta)$-MMP. Then there is a $(K_{\Ff'}+\Delta')$-MMP. In particular, there is a $(K_\Ff+\Delta)$-MMP.
\end{lem}

\begin{proof} By Lemma \ref{lem: contract reduce sing of f}, $(X',\Ff',\Delta')\in\mathcal{D}$. Possibly replacing $X$, $\Ff$ and $\Delta$ with $X'$, $\Ff'$ and $\Delta'$ respectively, we may assume that $\rho=\id_X$. We construct a $(K_{\Ff}+\Delta)$-MMP in the following way:
\begin{itemize}
    \item[\textbf{Step 1}] If $K_{\Ff}+\Delta$ is nef, there is nothing left to prove.
    \item[\textbf{Step 2}]  Suppose that $K_{\Ff}+\Delta$ is not nef, then there exists a $(K_{\Ff}+\Delta)$-negative extremal ray $C$.  By Lemma \ref{lem: contraction of negative extremal ray foliation}, Lemma \ref{lem: property of negative extremal ray foliated pairs}(1)(2)(4) and Lemma \ref{lem: contract reduce sing of f}, there is a contraction $f: (X,\Ff,\Delta)\rightarrow (X'',\Ff'',\Delta'')$ of $C$, such that $(X'',\Ff'',\Delta'')\in\mathcal{D}$.
    \item[\textbf{Step 3}] Replace $(X,\Ff,\Delta)$ with $(X'',\Ff'',\Delta'')$ and return to \textbf{Step 1}.
\end{itemize}
By Lemma \ref{lem: contract reduce sing of f}, the process \textbf{Step 1}-\textbf{Step 3} terminates as the number of closed points belonging to $\Sing(\Ff'')$ is strictly less than the number of closed points belonging to $\Sing(\Ff)$. Therefore we get a $(K_\Ff+\Delta)$-MMP.
\end{proof}

To further study the structure of the minimal model program for foliated surface pairs, we need to study the structure of the curves contracted by the minimal model program. It turns out that artificial chains and non-artificial chains behave in different ways. More precisely, for non-artificial chains, we have the following:

\begin{lem}\label{lem: structure of exc no artificial chain} Let $(X,\Ff,\Delta)\in\mathcal{D}$ be a foliated triple and $\rho: X\rightarrow X'$ a $(K_\Ff+\Delta)$-MMP, such that 
\begin{itemize}
    \item  $X$ does not contain any $(K_\Ff+\Delta)$-artificial chain, and
    \item  $(X',\Ff':=\rho_*\Ff,\Delta':=\rho_*\Delta)\in\mathcal{D}$.
\end{itemize} 
Then $\Ex(\rho)$ is the set of all the maximal $(K_\Ff+\Delta)$-chains. 
\end{lem}

\begin{proof}

\noindent\textbf{Step 1}.
We may assume that $K_{\Ff}+\Delta$ is not nef, otherwise there is nothing left to prove. Since $\rho$ is a $(K_{\Ff}+\Delta)$-MMP, by Lemma \ref{lem: contract reduce sing of f}, there is an integer $n>0$, a sequence of birational contractions
$$(X_0,\Ff_0,\Delta_0)\xrightarrow{\phi_1} (X_1,\Ff_1,\Delta_1)\xrightarrow{\phi_2}\dots\xrightarrow{\phi_n}(X_n,\Ff_n,\Delta_n),$$
irreducible curves $C_1,\dots,C_n$ on $X_0$, and an irreducible curve $C_i'$ on $X_{i-1}$ for every $1\leq i\leq n$ satisfying the following.
\begin{itemize}
    \item $(X,\Ff,\Delta)=(X_0,\Ff_0,\Delta_0)$ and $(X',\Ff',\Delta')=(X_n,\Ff_n,\Delta_n)$,
        \item $\rho=\phi_n\circ\dots\circ\phi_1$, and
        \item for every $1\leq i\leq n$,
        \begin{itemize}
            \item  $(X_i,\Ff_i,\Delta_i)\in\mathcal{D}$,
    \item $\phi_i$ is the contraction of $C_i'$, 
    \item $C_i'$ is a $(K_{\Ff_i}+\Delta_i)$-negative extremal ray, and
    \item $C_i'$ is the strict transform of $C_i$ on $X_{i-1}$.
        \end{itemize}
\end{itemize}
Let $V$ be the union of all the maximal $(K_\Ff+\Delta)$-chains. Since $\Ex(\rho)=\cup_{i=1}^nC_i$, in the rest of the proof, we only need to show that $\cup_{i=1}^nC_i=\Supp V$. \vspace{2mm}

\noindent\textbf{Step 2}. Since $\Supp V\subset\Ex(\rho)$, we only need to show that $\cup_{i=1}^nC_i\subset\Supp V$. We define
$$W:=\{i| 1\leq i\leq n, C_i\subset\Supp V\}.$$
If $W=\{1,2,\dots,n\}$ there is nothing left to prove. Otherwise, we may define
$$j:=\min\{i|1\leq i\leq n, C_i\not\in W\}.$$
By Lemma \ref{lem: property of negative extremal ray foliated pairs}(1), $C_j'$ is an $\Ff_{j-1}$-invariant curve. Thus $C_j$ is an $\Ff$-invariant curve.\vspace{2mm}

\noindent\textbf{Step 3}. We show that $Z(\Ff, C_j)=2$. By Lemma \ref{lem: property of negative extremal ray foliated pairs}(3), $Z(\Ff_{j-1},C_j')=1$. By Lemma \ref{lem: F chain and KF chain}(3), $Z(\Ff, C_j)\geq 1$. 

If $Z(\Ff, C_j)=1$, by Lemma \ref{lem: property of negative extremal ray foliated pairs}(2)(4)(5), $C_j$ itself is a $(K_\Ff+\Delta)$-chain, which contradicts to our assumptions.

If $Z(\Ff, C_j)\geq 3$, since there does not exist any $(K_{\Ff}+\Delta)$-artificial chain, the contraction of $C_i'$ is a cyclic quotient singularity on $X_i$ of order $\geq 2$ for every $1\leq i<j$.  Since $Z(\Ff_{j-1},C_j')=1$, $C_j'$ passes through at least $2$ cyclic quotient singularities of $X_{j-1}$, which contradicts to Lemma \ref{lem: property of negative extremal ray foliated pairs}(5). 

Thus $Z(\Ff, C_j)=2$.\vspace{2mm}

\noindent\textbf{Step 4}. Since $Z(\Ff, C_j)=2$ and $Z(\Ff_{j-1},C_j')=1$, there exists $1\leq k<j$ such that $C_k\cap C_j\not=\emptyset$. By Lemma \ref{lem: property of negative extremal ray foliated pairs}(1), $p:=C_k\cap C_j$ is a reduced singularity of $\Ff$. 

By Lemma \ref{lem: property of negative extremal ray foliated pairs}(5), $C_j'$ passes through at most $1$ closed point on $\Sing(X_{j-1})$. Since $X$ does not contain any $(K_\Ff+\Delta)$-artificial chain, for any integer $1\leq i<j$ such that $i\not=k$,  $C_i\cap C_j=\emptyset$.

Since $C_k$ belongs to a $(K_\Ff+\Delta)$-chain and $C_k$ and $C_j$ are the only two $\Ff$-invariant curves passing through $p$, there exists an integer $1\leq l<j$ and a  $(K_\Ff+\Delta)$-chain 
$$D:=\cup_{i=1}^{l}C_{k_i}$$
where $C_{k_l}=C_k$ and $k_i<k_{i+1}$ for any $1\leq i\leq l-1$. By applying Lemma \ref{lem: Intersection surface} and Lemma \ref{lem: property of negative extremal ray foliated pairs}(1)(5), we deduce that $D\cup C_j$ is a $(K_{\Ff}+\Delta)$-chain, which contradicts to $C_j\not\in W$.
\end{proof}

The next lemma gives a precise description on the structure of non-artificial chains and artificial chains of the minimal model program of foliated pairs. Generally speaking, a minimal model program of foliated pairs can be visualized as the composition of a partial minimal model program which contracts exactly all the non-artificial chains, and a partial minimal model program which only contracts artificial chains.

\begin{lem}\label{lem: mmp foliated pairs} Let $(X:=X_0,\Ff:=\Ff_0,\Delta:=\Delta_0)\in\mathcal{D}$ be a foliated triple. Then there is 
\begin{itemize}
    \item an integer $n\geq 0$,
    \item a $(K_\Ff+\Delta)$-MMP $\rho: X\rightarrow X'$,
    \item a sequence of birational contractions
$$(X_0,\Ff_0,\Delta_0)\xrightarrow{\phi_1} (X_1,\Ff_1,\Delta_1)\xrightarrow{\phi_2}\dots\xrightarrow{\phi_n}(X_n,\Ff_n,\Delta_n)$$
and
    \item a contraction $\nu: X_n\rightarrow X'$
\end{itemize} 
satisfying the following. 
\begin{enumerate}
\item $(X',\Ff':=\rho_*\Ff,\Delta':=\rho_*\Delta)\in\mathcal{D}$,
\item for every $1\leq i\leq n$, $(X_i,\Ff_i,\Delta_i)\in\mathcal{D}$, 
\item $\rho=\nu\circ\phi_n\circ\dots\circ\phi_1$,
\item for every $1\leq i\leq n$, $\phi_i$ is a contraction of $(K_{\Ff_{i-1}}+\Delta_{i-1})$-artificial chains,
\item $X_n$ does not contain any $(K_{\Ff_n}+\Delta_n)$-artificial chain, and
\item $\nu$ exactly contracts all the maximal $(K_{\Ff_n}+\Delta_n)$-chains.
\end{enumerate}

\end{lem}

\begin{proof}
First we find $n$ which satisfies (2)(4)(5) in the following way.\vspace{2mm}

\noindent\textbf{Step 1}. Suppose that we have already defined $X_{i-1},\Ff_{i-1}$ and $\Delta_{i-1}$ for some integer $i>0$. If $X_{i-1}$ does not contain any $(K_{\Ff_{i-1}}+\Delta_{i-1})$-artificial chain, we let $n:=i-1$ and the lemma follows from Lemma \ref{lem: structure of exc no artificial chain}.\vspace{2mm}

\noindent\textbf{Step 2}. If $X_{i-1}$ contains some $(K_{\Ff_{i-1}}+\Delta_{i-1})$-artificial chains, then we let $\phi_i: (X_{i-1},\Ff_{i-1},\Delta_{i-1})\rightarrow (X_{i},\Ff_{i},\Delta_{i})$ be the contraction of all the  $(K_{\Ff_{i-1}}+\Delta_{i-1})$-artificial chains.\vspace{2mm}

\noindent\textbf{Step 3}. (2)(4) hold immediately. Replacing $i$ with $i+1$, we repeat this process from \textbf{Step 1}.\vspace{2mm}

By Lemma \ref{lem: contract reduce sing of f}, the process \textbf{Step 1}-\textbf{Step 3} terminates as the number of closed points belonging to $\Sing(\Ff_i)$ is strictly less than the number of closed points belonging to $\Sing(\Ff_{i-1})$ for any $i$.

By Lemma \ref{lem: mmp foliated pairs weak version}, there exists a $(K_{\Ff_n}+\Delta_n)$-MMP $\nu:(X_n,\Ff_n,\Delta_n)\rightarrow (X',\Ff',\Delta')$. We let $\rho:=\nu\circ\phi_n\circ\dots\circ\phi_1$, and the lemma follows from Lemma \ref{lem: contract reduce sing of f} and Lemma \ref{lem: structure of exc no artificial chain}.
 \end{proof}

\section{Zariski decomposition of foliated surface pairs} In this section, we study the Zariski decomposition of foliated surface pairs, which is the main ingredient in the proof of the invariance of plurigenera (Theorem \ref{thm: foliated surface pair inv plur}). Moreover, we deduce a vanishing theorem at the end of this subsection, which immediately implies Theorem \ref{thm: weak kv vanishing foliated surface}.

\begin{lem}\label{lem: negative part under birational contraction} Let $(X,\Ff,\Delta)\in\mathcal{D}$ be a foliated triple and $\rho: (X,\Ff,\Delta)\rightarrow (X',\Ff',\Delta')$ a $(K_{\Ff}+\Delta)$-MMP. Assume that
$$K_\Ff+\Delta=\rho^*(K_{\Ff'}+\Delta')+N,$$
and let $P:=\rho^*(K_{\Ff'}+\Delta')$. Then
\begin{enumerate}
\item $P$ is nef,
\item $P\cdot C=0$ for any curve $C\subset\Supp N$,
\item $N\geq 0$, and
\item $\Supp N=\Supp\Ex(\rho)$.
\end{enumerate}
\end{lem}

\begin{proof} Since $K_\Ff+\Delta$ is pseudo-effective, $K_{\Ff'}+\Delta'$ is nef, which implies (1). Since $P\cdot C=0$ for any curve $C\subset\Supp\Ex(\rho)$ and $\Supp N\subset\Supp\Ex(\rho)$, we deduce (2). Since $\rho$ is a $(K_{\Ff}+\Delta)$-MMP, $\rho$ is $(K_{\Ff}+\Delta)$-negative, and (3)(4) are immediate.
\end{proof}

The next lemma shows that for any foliation without any artificial chain, $P$ and $N$ as in Lemma \ref{lem: negative part under birational contraction} are exactly the positive part and the negative part of the Zariski decomposition of $K_{\Ff}+\Delta$.

\begin{lem}\label{lem: structure of zdecomposition of foliation} Assume that
\begin{itemize}
    \item $(X,\Ff,\Delta)\in\mathcal{D}$ is a foliated triple such that $X$ does not contain any $(K_{\Ff}+\Delta)$-artificial chain,
    \item $\rho: (X,\Ff,\Delta)\rightarrow (X',\Ff',\Delta')$ is a $(K_{\Ff}+\Delta)$-MMP, and
    \item $P:=\rho^*(K_{\Ff'}+\Delta')$ and $N:=(K_\Ff+\Delta)-P$,
\end{itemize}
then
\begin{enumerate}
\item $\Supp N$ is the union of all the maximal $(K_\Ff+\Delta)$-chains, 
\item each $(K_\Ff+\Delta)$-chain has a negative definite intersection matrix, 
\item the intersection matrix of $N$ is negative definite, and
\item $P$ and $N$ are the positive part and the negative part of the Zariski decomposition of $K_\Ff+\Delta$ respectively.
\end{enumerate}
\end{lem}

\begin{proof} We may assume that there exists a $(K_\Ff+\Delta)$-chain, otherwise $X=X'$ and $N=0$ and there is nothing left to prove.

(1) follows from Lemma \ref{lem: structure of exc no artificial chain} and Lemma \ref{lem: negative part under birational contraction}(4). 

To prove (2), let $n>0$ be an integer and $C=\cup_{i=1}^nC_i$  a $(K_{\Ff}+\Delta)$-chain. Then there is a sequence of contractions
$$(X_0,\Ff_0,\Delta_0):=(X,\Ff,\Delta)\xrightarrow{\phi_1} (X_1,\Ff_1,\Delta_1)\xrightarrow{\phi_2}\dots\xrightarrow{\phi_n}(X_n,\Ff_n,\Delta_n),$$
and an irreducible curve $C_i'$ on $X_{i-1}$ for every $1\leq i\leq n$ satisfying the following.
\begin{itemize}
    \item $\phi_i$ is the contraction of $C_i'$,
    \item $C_i'$ is the birational transform of $C_i$ on $X_{i-1}$, and
    \item $(K_{\Ff_{i-1}}+\Delta_{i-1})\cdot C_i'<0$.
\end{itemize}
Since $K_{\Ff}+\Delta$ is pseudo-effective, $K_{\Ff_{i-1}}+\Delta_{i-1}$ is pseudo-effective. Thus $(C_{i}')^2<0$. 
By Lemma \ref{lem: Intersection surface}, for any $1\leq i\leq n-1$, 
$(C_{i+1}')^2=C_{i+1}^2-\frac{1}{(C_i')^2}.$
Thus for any non-zero real vector $\bm{x}=(x_1,..., x_n)$, by a complicated but elementary calculation, we have
$$\bm{x}||C||\bm{x}^T=-(\sqrt{-(C_{n}')^2}x_n)^2-\sum_{i=1}^{n-1}(\sqrt{-(C_{i}')^2}x_i+\frac{1}{\sqrt{-(C_{i}')^2}}x_{i+1})^2<0.$$
Thus $||C||$ is negative definite, and we deduce (2). (3) follows from (2). (4) follows from (3) and Lemma \ref{lem: negative part under birational contraction}.\end{proof}


The next lemma shows that the Zariski decomposition of minimal foliated pairs is ``stable" under small perturbation by ample divisors. More precisely, we have the following:

\begin{lem}\label{lem: zd of perturbation with small ample} Assume that
\begin{itemize}
    \item $(X,\Ff,\Delta)\in\mathcal{D}$ is a foliated triple such that $X$ does not contain any $(K_{\Ff}+\Delta)$-artificial chain,
    \item $P$ and $N$ are the positive part and the negative part of the Zariski decomposition of $K_{\Ff}+\Delta$ respectively,
    \item $A$ is an ample $\Qq$-divisor on $X$,
    \item for every rational number $\epsilon\geq 0$,  $P_{\epsilon}^A$ and $N^A_{\epsilon}$ are the positive part and the negative part of the Zariski decomposition of $K_{\Ff}+\Delta+\epsilon A$ respectively,
\end{itemize}
then there is a rational number $\epsilon_0>0$ satisfying the following. For any real number $0\leq\epsilon\leq\epsilon_0$,
\begin{enumerate}
\item $\Supp N^A_{\epsilon}=\Supp N$,
 \item any $(K_\Ff+\Delta)$-MMP is also a $(K_\Ff+\Delta+\epsilon A)$-MMP, and
 \item $N\geq N^A_{\epsilon}$.
\end{enumerate}
\end{lem}

\begin{proof} Possibly replacing $A$ with a general element of $|A|_{\Qq}$, we may assume that $A\geq 0$. By Lemma \ref{lem: structure of exc no artificial chain}, any $(K_{\Ff}+\Delta)$-MMP is a contraction of all the maximal $(K_{\Ff}+\Delta)$-chains. Thus there are only finitely many foliated triples constructed in all the possible partial $(K_{\Ff}+\Delta)$-MMPs, and we may suppose that they are $\{(X_i,\Ff_i,\Delta_i)\}_{i=1}^n$. By Lemma \ref{lem: contract reduce sing of f}, each $(X_i,\Ff_i,\Delta_i)\in\mathcal{D}$.

Let $A_i$ be the birational transform of $A$ on $X_i$ for each $i$. By Lemma \ref{lem: Intersection surface}, $A_i$ is ample. Therefore,
\begin{itemize}
    \item if $K_{\Ff_i}+\Delta_i$ is nef, then $(K_{
\Ff_i}+\Delta_i+\epsilon A_i)$ is also nef for any $\epsilon>0$, and
\item there exists $0<\epsilon\ll 1$, such that any $(K_{\Ff_i}+\Delta_i)$-negative extremal ray is also a $(K_{\Ff_i}+\Delta_i+\epsilon A_i)$-negative extremal ray for every $i$. 
\end{itemize}
Thus we deduce (2). (1) follows from (2) and Lemma \ref{lem: negative part under birational contraction}.

For any irreducible component $C$ of $\Supp N^A_{\epsilon}=\Supp N$, since $P^{A}_{\epsilon}\cdot C=P\cdot C=0$, we have
\begin{align*}
(N^A_{\epsilon}-N)\cdot C&=((P^A_{\epsilon}+N^A_{\epsilon}-(P+N))+(P-P^A_{\epsilon}))\cdot C\\
&=(P^A_{\epsilon}+N^A_{\epsilon}-(P+N))\cdot C\\
&=(K_{\Ff}+\Delta+\epsilon A-(K_{\Ff}+\Delta))\cdot C\\
&=\epsilon A\cdot C>0.
\end{align*}
Let $X\rightarrow X'$ be the contraction of $\Supp N^A_{\epsilon}=\Supp N$. Then $N^A_{\epsilon}-N$ is nef over $X'$. By the negativity Lemma \ref{lem: num negativity lemma}(2), $N\geq N^A_{\epsilon}$. \end{proof}

The next lemma shows that under a small perturbation, the null locus of the positive part of the Zariski decomposition is exactly the support of the negative part of the Zariski decomposition for foliated surface pairs:

\begin{lem}\label{lem: structure of nullA} Assume that
\begin{itemize}
    \item $(X,\Ff,\Delta)\in\mathcal{D}$ is a foliated triple such that $X$ does not contain any $(K_{\Ff}+\Delta)$-artificial chain,
\item $P$ and $N$ are the positive part and the negative part of the Zariski decomposition of $K_{\Ff}+\Delta$ respectively,
    \item $A$ is an ample $\Qq$-divisor on $X$,
\item $0<\epsilon\ll 1$ is a sufficiently small rational number, and
\item  $P^A_{\epsilon}$ and $N^{A}_{\epsilon}$ are the positive part and the negative part of the Zariski decomposition of $K_{\Ff}+\Delta+\epsilon A$ respectively.
\end{itemize}
Then $\Null(P^A_{\epsilon})=\Supp N.$
\end{lem}
\begin{proof}
Let $\rho: X\rightarrow X'$ be a $(K_{\Ff}+\Delta)$-MMP. By  Lemma \ref{lem: zd of perturbation with small ample}(2), $\rho$ is a $(K_{\Ff}+\Delta+\epsilon A)$-MMP. By Lemma \ref{lem: Intersection surface}, $\rho_*A$ is ample. For any curve $C\subset \Null(P^A_{\epsilon})$, we have 
\begin{align*}
0&\leq\epsilon\rho_*A\cdot\rho_*C\leq\rho_*(K_{\Ff}+\Delta)\cdot\rho_*C+\epsilon\rho_* A\cdot\rho_*C\\
&=\rho_*(K_{\Ff}+\Delta+\epsilon A)\cdot\rho_*C=\rho^*\rho_*(K_{\Ff}+\Delta+\epsilon A)\cdot C\\
&=P^A_{\epsilon}\cdot C=0.
\end{align*}
Thus $\rho_*A\cdot\rho_*C=0$, hence $C$ is $\rho$-exceptional. By Lemma \ref{lem: negative part under birational contraction}, $C\subset\Supp N$. Thus $\Null(P^A_{\epsilon})\subset\Supp N$. By Lemma \ref{lem: zd of perturbation with small ample}(1), $\Supp N=\Supp N^A_{\epsilon}\subset \Null(P^A_{\epsilon})$. Thus $\Null(P^A_{\epsilon})=\Supp N.$
\end{proof}

Now we turn to study the structure of the null locus of the positive part of the Zariski decomposition of foliated pairs without the perturbation by an ample divisor. To simplify our statements, we may only consider foliated triples $(X,\Ff,\Delta)\in\mathcal{D}_0$.

\begin{lem}\label{lem: structure of null} Assume that
\begin{itemize}
    \item $(X,\Ff,\Delta)\in\mathcal{D}_0$ is a foliated triple, and
    \item $P$ and $N$ are the positive part and the negative part of the Zariski decomposition of $K_{\Ff}+\Delta$ respectively,
\end{itemize}
then for any irreducible curve $C\subset\Null(P)$, 
\begin{enumerate}
    \item $C$ is $\Ff$-invariant,
    \item  $C^2<0$,
\end{enumerate} 
and $C$ is of one of the following types:
\begin{itemize}
\item[\textbf{\rm \textbf{Type A}}] $C\in\Supp N$.
\item[\textbf{\rm \textbf{Type B}}] $C$ is a nodal curve, whose unique singularity is a reduced singularity of $\Ff$. Moreover, 
\begin{itemize}
\item[(B.1)] $C$ is a connected component of $\Null(P)$, and
\item[(B.2)] $K_\Ff\cdot C=\Delta\cdot C=0$.
\end{itemize}
\item[\textbf{\rm \textbf{Type C}}] $C$ is a component of a cycle of rational curves, i.e. there exists an integer $n\geq 2$, smooth rational curves $C_1:=C,C_2\dots,C_n$, such that 
\begin{itemize}
\item[(C.1)] $\cup_{i=1}^{n}C_i$ is a connected component of $\Null(P)$, 
\item[(C.2)] for any integers $1\leq i,j\leq n$, 
\begin{itemize}
\item $C_i\cdot C_j=1$ if $|i-j|=1$ or $\{i,j\}=\{1,n\}$, and
\item  $C_i\cdot C_j=0$ otherwise,
\end{itemize}
and
\item[(C.3)]  $K_\Ff\cdot C=\Delta\cdot C=0$.
\end{itemize}
\item[\textbf{\rm \textbf{Type D}}] $C$ is a smooth rational curve, such that
\begin{itemize}
\item[(D.1)] $C\cap\Supp N=\emptyset$, 
\item[(D.2)] $Z(\Ff,C)=1$ or $2$,
\item[(D.3)] $C$ is not of \textbf{\rm \textbf{Type C}}.
\end{itemize}
\item[\textbf{\rm \textbf{Type E}}] $C$ is a smooth rational curve, such that
\begin{itemize}
\item[(E.1)] $C$ is the tail of a $(K_\Ff+\Delta)$-chain $D$ such that $D\cup C$ is an $\Ff$-chain, and
\item[(E.2)] $C$ does not intersect any $(K_\Ff+\Delta)$-chain except $D$.
\end{itemize}
\item[\textbf{\rm \textbf{Type F}}] $C$ is a smooth rational curve, such that
\begin{itemize}
\item[(F.1)] $C$ is the tail of two $(K_\Ff+\Delta)$-chains $D_1$ and $D_2$, such that $\det(-||D_1||)=\det(-||D_2||)=2$, 
\item[(F.2)] $C$ does not intersect any $(K_\Ff+\Delta)$-chain except $D_1$ and $D_2$, and
\item[(F.3)] $C$ intersects at most one other irreducible $\Ff$-invariant curve which does not belong to $D_1$ and $D_2$. 
\end{itemize}
\end{itemize}
Moreover, for any curve $C$ of \textbf{\rm \textbf{Type D}} or \textbf{\rm \textbf{Type E}} or \textbf{\rm \textbf{Type F}}, the connected component of $\Null(P)$ containing $C$ is a string.
\end{lem}

\begin{proof} Let $C\subset\Null(P)$ be an irreducible curve. If $C\subset\Supp N$, then $C$ is of \textbf{Type A}. By Lemma \ref{lem: structure of exc no artificial chain} and Lemma \ref{lem: negative part under birational contraction}, (1)(2) hold in this case. Thus in the rest of the proof, we may assume that $C\not\subset\Supp N$.

Let $\rho: X\rightarrow X'$ be a $(K_{\Ff}+\Delta)$-MMP, $\Ff':=\rho_*\Ff$, $\Delta':=\rho_*\Delta$, and $C':=\rho_*C$. 

First we show (1)(2) and $\chi(C')\geq Z(\Ff',C')$. Since $C\not\subset\Supp N$, by Lemma \ref{lem: negative part under birational contraction}, $C'\not=0$. Since $P\cdot C=0$ and $P=\rho^*(K_{\Ff'}+\Delta')$, we have $(K_{\Ff'}+\Delta')\cdot C'=0$. Since $\Delta\wedge\textbf{\rm \textbf{B}}_{+}(K_\Ff+\Delta)=0$, by Lemma \ref{lem: aug stable base locus pushforward},  $\Delta'\wedge\textbf{\rm \textbf{B}}_{+}(K_{\Ff'}+\Delta')=0$. Thus $\Delta'\cdot C'\geq 0$. Since $K_{\Ff}+\Delta$ is big, $K_{\Ff'}+\Delta'$ is big. Thus $(C')^2<0$. We deduce that
$$(K_{\Ff'}+C')\cdot C'<K_{\Ff'}\cdot C'\leq (K_{\Ff'}+\Delta')\cdot C'=0.$$
By Lemma \ref{lem: tang formula}, $C'$ is $\Ff'$-invariant. Thus $C$ is $\Ff$-invariant. Since $(C')^2<0$, by Lemma \ref{lem: Intersection surface}, $C^2<0$. By Lemma \ref{lem: z and cs formula}(1),
$$\chi(C')=Z(\Ff',C')-K_{\Ff'}\cdot C'\geq Z(\Ff',C').$$ 

 To classify the types of $C$, there are two cases:
 \begin{itemize}
     \item[\textbf{Case 1}] $C\cap\Supp N\not=\emptyset$, and
     \item[\textbf{Case 2}] $C\cap\Supp N=\emptyset$.
 \end{itemize}
 
 \noindent\textbf{Case 1}. In this case, we show that $C$ is of \textbf{Type E} or \textbf{Type F}. 
 
  Since $X$ does not contain any $(K_{\Ff}+\Delta)$-artificial chain and $C\cap\Supp N\not=\emptyset$, $C'$ contains at least $1$ cyclic quotient singularity of $X'$ of order $\geq 2$. Since $Z(\Ff',C')\geq 0$, $\chi(C')\geq 0$. Thus $C'$ is a rational curve, which implies that $C$ is a smooth rational curve.

 Since $(C')^2<0$, by Lemma \ref{lem: z and cs formula}(2), $CS(\Ff',C')<0$. By Lemma \ref{lem: basic properties of z and cs}(3), $C$ passes through at least $1$ singularity of $\Ff'$. By Lemma \ref{lem: basic properties of z and cs}(4), $Z(\Ff',C')\geq 1$. Thus $\chi(C')\geq 1$.
 
 Let $m\geq 0$ be an integer such that $C$ is the tail of exactly $m$ $(K_\Ff+\Delta)$-chains. Since $X$ does not contain any $(K_{\Ff}+\Delta)$-artificial chain and since $\rho$ contracts all the $(K_\Ff+\Delta)$-chains, we have
$$\chi(C')=2-\sum_{i=1}^{m}(1-\frac{1}{r_i})$$
where the $r_i$ are the orders of cyclic quotient singularities of $X'$ along $C'$ such that each $r_i\geq 2$.  Since $C\cap\Supp N\not=\emptyset$, $m>0$. Thus $\chi(C')<2$, which implies that $Z(\Ff',C')=1$. Therefore $\chi(C')\geq 1$, and there are two cases:

\begin{itemize}
    \item[\textbf{Case 1.1}] $m=1$.
    \item[\textbf{Case 1.2}] $m=2$ and $r_1=r_2=2$.
\end{itemize}

\noindent\textbf{Case 1.1} In this case, by Lemma \ref{lem: F chain and KF chain}(5), $Z(\Ff,C)=2$, and we deduce that $C$ is of \textbf{Type E}.\vspace{2mm}

\noindent\textbf{Case 1.2} In this case, $\chi(C')=1$. Let $D_1$ and $D_2$ be the $(K_\Ff+\Delta)$-chains which intersect $C$. By Lemma \ref{lem: z and cs formula}(1) and Lemma \ref{lem: contraction intersection numbers},
\begin{align*}
0&=-\chi(C')+Z(\Ff',C')=K_{\Ff'}\cdot C'\\
&=K_{\Ff}\cdot C-\det(-||D_1||)^{-1}-\det(-||D_2||)^{-1}.
\end{align*}
Sine $X$ is smooth, $K_{\Ff}\cdot C\in\mathbb Z$ and $\det(-||D_1||),\det(-||D_2||)\in\mathbb N^+$. Since $D_1$ and $D_2$ are not artificial $\Ff$-chains, $\det(-||D_1||),\det(-||D_2||)\geq 2$. Thus $\det(-||D_1||)^{-1}=\det(-||D_2||)^{-1}=\frac{1}{2}$. Since $Z(\Ff',C')=1$, $C'$ intersects at most one other $\Ff'$-invariant curve, and we deduce that $C$ is of \textbf{Type F}.\vspace{2mm}

\noindent\textbf{Case 2}. In this case, $\rho$ is an isomorphism in an open neighborhood of $C$. In particular,
\begin{itemize}
\item $(K_\Ff+\Delta)\cdot C=(K_{\Ff'}+\Delta')\cdot C'=0$, and
\item $\chi(C)=\chi(C')\geq Z(\Ff',C')=Z(\Ff,C)\geq 0$.
\end{itemize}

Since $X$ is smooth, there are two cases:
\begin{itemize}
    \item[\textbf{Case 2.1}] $\chi(C)=0$.
    \item[\textbf{Case 2.2}] $\chi(C)=2$.
\end{itemize}

\noindent\textbf{Case 2.1}. In this case, $Z(\Ff,C)=0$. By Lemma \ref{lem: z and cs formula}(1), $K_\Ff\cdot C=0$, hence $\Delta\cdot C=0$. By Lemma \ref{lem: z and cs formula}(2), since $C^2<0$, $CS(\Ff,C)<0$. By Lemma \ref{lem: basic properties of z and cs}(3), $C\cap \Sing(\Ff)\not=\emptyset$. By Lemma \ref{lem: basic properties of z and cs}(4), $C$ is not smooth. Thus $C$ is a nodal curve, whose singularity $p$ is a singularity of $\Ff$. Since $\Ff$ is canonical, by Lemma \ref{lem: eigenvalue and cs and z formula}(4), $p$ is not a Poincar\'{e}-Dulac singularity. By Lemma \ref{lem: either type b pd or reduced}, $p$ is a reduced singularity.

Since $Z(\Ff,C)=Z(\Ff,C,p)=0$, by Lemma \ref{lem: basic properties of z and cs}(4), $p=C\cap\Sing(\Ff).$ In particular, $C$ does not intersect any other $\Ff$-invariant curve. Therefore $C$ is of \textbf{Type B}.\vspace{2mm}

\noindent\textbf{Case 2.2}. If $\chi(C)=2$, then $C$ is a smooth rational curve, and 
$$0\leq Z(\Ff,C)=K_{\Ff}\cdot C+\chi(C)\leq (K_{\Ff}+\Delta)\cdot C+2=2.$$
Thus $Z(\Ff,C)=1$ or $2$, which implies that $C$ is of \textbf{Type C} or \textbf{Type D}.

For any curve $C$ of \textbf{Type D} or \textbf{Type E} or \textbf{Type F}, let $G$ be the connected component of $\Null(P)\backslash\Supp N$ containing $C$. Then either $G\cap\Supp N=\emptyset$ or $G\cap\Supp N\not=\emptyset$.
\begin{itemize}
\item 
Suppose that $G\cap\Supp N=\emptyset$. Since $C$ is not of \textbf{Type C}, each irreducible component $C''$ of $G$ is of \textbf{Type D}. In particular $Z(\Ff,C'')=1$ or $2$. By Lemma \ref{lem: basic properties of z and cs}(4), $C''$ intersects at most $2$ irreducible components of $G\backslash C''$. In particular, $G$ is either a string or a cycle. If $G$ is a cycle, then $K_{\Ff}\cdot C=-\chi(C)+Z(\Ff,C)=0$ by Lemma \ref{lem: z and cs formula}(1), which implies that $C$ is of \textbf{Type C}, which is not possible. 
\item Suppose that $G\cap\Supp N\not=\emptyset$. Then each irreducible component $C''$ of $G$ is of \textbf{Type D}, \textbf{Type E} or \textbf{Type F}. We have the following:
\begin{itemize}
    \item If $C''$ is of \textbf{Type E} or \textbf{Type F}, by (E.1) and (F.3), $C''$ intersects at most $1$ irreducible component of $G\backslash C''$. 
    \item If $C''$ is of \textbf{Type D}, then $Z(\Ff,C'')=1$ or $2$. By Lemma \ref{lem: basic properties of z and cs}(4), $C''$ intersects at most $2$ irreducible components of $G\backslash C''$.
\end{itemize}
\end{itemize}
Thus $G$ is a string.
\end{proof}

We are now ready to give a vanishing theorem for foliated surface pairs. Before we state the theorem, we first prove a technical lemma which provides an auxiliary pair in the proof of the vanishing theorem:

\begin{lem}\label{lem: find pair for kv vanishing} Assume that
\begin{itemize}
    \item $(X,\Ff,\Delta)\in\mathcal{D}_0$ is a foliated triple, and
    \item $P$ and $N$ are the positive part and the negative part of the Zariski decomposition of $K_{\Ff}+\Delta$ respectively,
\end{itemize}
then there exists a $\mathbb Q$-divisor $\Theta$ on $X$, such that 
\begin{enumerate}
\item $\Supp\Theta\subset\Supp N$,
\item $\lfloor\Theta\rfloor=0$,
\item for any maximal $(K_\Ff+\Delta)$-chain $C=\cup_{i=1}^{n}C_i$, 
\begin{itemize}
\item for any $1\leq i\leq n-1$, $(K_X+\Theta)\cdot C_i\leq 0$, and
\item $(K_X+\Theta)\cdot C_n\leq -1$, 
\end{itemize}
and
\item for any irreducible curve $E$ such that $E^2=-2$ and $E$ is a maximal $(K_\Ff+\Delta)$-chain, $\mu_E\Theta=\frac{1}{2}$.
\end{enumerate}
 \end{lem}

\begin{proof} For any $1\leq i\leq n$, let $c_i:=\min\{2+C_i^2,0\}.$  Since $X$ is smooth and $C_i$ is a smooth rational curve, 
$$K_X\cdot C_i=-2-C_i^2.$$

Let $\lambda_0:=0$ and $\lambda_{n+1}=1$. We consider the solutions $(\lambda_1,\dots,\lambda_n)$ of the equations
\begin{equation}\label{equ: find pair for kv vanishing}
c_i=\lambda_{i-1}+(c_i-2)\lambda_i+\lambda_{i+1}, \forall 1\leq i\leq n.\tag{5.6.1}
\end{equation}

\noindent\textbf{Fact 1}. (5.6.1) has a solution.
\begin{proof}[Proof of \textbf{Fact 1}]
We define $\bm{A}:=\{c_{i,j}\}_{1\leq i,j\leq n}$ to be the following $n\times n$ matrix:
\begin{itemize}
\item For any $1\leq i\leq n$, $c_{i,i}=c_i-2$.
\item For any $1\leq i,j\leq n$ such that $|i-j|=1$, $c_{i,j}=1$.
\item For any $1\leq i,j\leq n$ such that $|i-j|\geq 2$, $c_{i,j}=0$.
\end{itemize}
Then (5.6.1) is equivalent to the equality
$$(c_1,c_2,\dots, c_{n-1},c_n-1)=(\lambda_1,\dots,\lambda_n)\cdot\bm{A}.$$
Since $c_{i,i}\leq -2$ for any $1\leq i\leq n$, $\bm{A}$ is negative definite. In particular, $\bm{A}$ is invertible, and we have
$$(\lambda_1,\dots,\lambda_n)=(c_1,c_2,\dots, c_{n-1},c_n-1)\cdot\bm{A}^{-1}.$$

\end{proof}

\noindent\textbf{Fact 2}. $0\leq\lambda_i<1$ for any $1\leq i\leq n$. 
\begin{proof}[Proof of \textbf{Fact 2}]
Since
$$(c_1,c_2,\dots, c_{n-1},c_n-1)=(\lambda_1,\dots,\lambda_n)\cdot\bm{A},$$
and $c_i\leq 0$ for any $1\leq i\leq n$, by Lemma \ref{lem: num negativity lemma}(1), $\lambda_i\geq 0$ for any $1\leq i\leq n$. Moreover, since
$$(c_1-1,c_2,\dots,c_{n-1},c_n-1)=(1,1,\dots,1)\cdot\bm{A},$$
we have
$$(-1,0,\dots,0)=(1-\lambda_1,\dots,1-\lambda_n)\bm{A}.$$
By Lemma \ref{lem: num negativity lemma}(1), $\lambda_i\leq 1$ for every $1\leq i\leq n$. 

If there exists $1\leq k\leq n$ such that $\lambda_k=1$, then 
$$c_k=\lambda_{k-1}+(c_k-2)+\lambda_{k+1},$$
which implies that
$$\lambda_{k-1}+\lambda_{k+1}=2.$$
Therefore $\lambda_{k-1}=\lambda_{k+1}=1$. We inductively deduce that $\lambda_i=1$ for any $1\leq i\leq n$, which contradicts to $c_1=(c_1-2)\lambda_1+\lambda_2$.
\end{proof}

\noindent\textbf{Fact 3}. If $n=1$ and $C_1^2=-2$, then $c_1=0$ and $\bm{A}=\{-2\}$, hence $\lambda_1=\frac{1}{2}.$\vspace{2mm}

We define
$$\Theta_C:=\sum_{i=1}^{n}\lambda_iC_i.$$
By the three facts above, we have
\begin{itemize}
\item $\Supp\Theta_C\subset C\subset\Supp N$,
\item $\lfloor\Theta_C\rfloor=0$, 
\item if $1\leq i\leq n-1$ and $C_i^2\leq -2$, then
$$\Theta_C\cdot C_i=\lambda_{i-1}+(C_i^2)\lambda_i +\lambda_{i+1}=\lambda_{i-1}+(c_i-2)\lambda_i+\lambda_{i+1}=c_i=2+C_i^2,$$
\item  if $1\leq i\leq n-1$ and $C_i^2=-1$, then
$$\Theta_C\cdot C_i=\lambda_{i-1}+(C_i^2)\lambda_i +\lambda_{i+1}=\lambda_i+(\lambda_{i-1}-2\lambda_i+\lambda_{i+1})=\lambda_i+c_i=\lambda_i,$$
\item if $C_n^2\leq -2$, then
$$\Theta_C\cdot C_n=\lambda_{n-1}+(C_n^2)\lambda_n=\lambda_{n-1}+(c_n-2)\lambda_n+(\lambda_{n+1}-1)=c_n-1=1+C_n^2,$$
\item if $C_n^2=-1$, then
$$\Theta_C\cdot C_n=\lambda_{n-1}+(C_n^2)\lambda_n=\lambda_{n-1}-2\lambda_n+\lambda_{n}+(\lambda_{n+1}-1)=c_n+\lambda_n-1=\lambda_n-1,$$
and
\item if $n=1$ and $C_1^2=-2$, then $\Theta_C=\frac{1}{2}C_1$.
\end{itemize}
Thus 
\begin{itemize}
\item if $1\leq i\leq n-1$ and $C_i^2\leq -2$, then
$$(K_X+\Theta_C)\cdot C=-2-C_i^2+2+C_i^2=0,$$
\item if $1\leq i\leq n-1$ and $C_i^2=-1$, then
$$(K_X+\Theta_C)\cdot C=\lambda_i-1<0,$$
\item if $C_n^2\leq -2$, then
$$(K_X+\Theta_C)\cdot C=1+C_n^2-(2+C_n^2)=-1,$$
and
\item if $C_n^2=-1$, then
$$(K_X+\Theta_C)\cdot C=\lambda_n-1-1<-1.$$
\end{itemize}

Therefore
$$\Theta:=\sum_{C\text{ is a maximal }(K_{\Ff}+\Delta)\text{-chain}}\Theta_C,$$
satisfies our requirements. \end{proof}

With this auxiliary divisor, we may prove the vanishing theorems. First we state a vanishing theorem of foliated surface pairs under the perturbation by some ample divisor:

\begin{lem}\label{lem: kv vanishing foliation version 1} Assume that 
\begin{itemize}
    \item $(X,\Ff,\Delta)\in\mathcal{D}$ is a foliated triple, such that
    \begin{itemize}
        \item $X$ is smooth, and
        \item $X$ does not contain any $(K_{\Ff}+\Delta)$-artificial chain,
    \end{itemize}
    \item $A$ is an ample $\Qq$-divisor on $X$,
    \item $0<\epsilon\ll 1$ is a sufficiently small rational number, and
    \item $P^A_{\epsilon}$ and $N^{A}_{\epsilon}$ are the positive part and the negative part of the Zariski decomposition of $K_{\Ff}+\Delta+\epsilon A$ respectively.
\end{itemize}
Then for any sufficiently divisible integer $m>0$ and any integer $i>0$,
$$H^i(X, mP^A_{\epsilon})=0.$$
\end{lem}

\begin{proof} Let $N$ be the negative part of the Zariski decomposition of $K_{\Ff}+\Delta$. By Lemma \ref{lem: zd of perturbation with small ample}(1), $\Supp N^A_{\epsilon}=\Supp N$. By Lemma \ref{lem: find pair for kv vanishing}, there exists an effective $\mathbb Q$-divisor $\Theta^A_{\epsilon}$ on $X$, such that
\begin{itemize}
\item $\Supp\Theta^A_{\epsilon}\subset\Supp N$,
\item $\lfloor\Theta^A_{\epsilon}\rfloor=0$, and
\item for any irreducible curve $C\subset\Supp N$, $(K_X+\Theta^A_{\epsilon})\cdot C\leq0$.
\end{itemize}

Since $A$ is ample and $K_\Ff+\Delta$ is pseudo-effective, $K_{\Ff}+\Delta+\epsilon A$ is big. Thus $P^A_{\epsilon}$ is big. In particular, there exists an integer $m_0>0$, such that $m_0P^A_{\epsilon}-(K_X+\Theta^A_{\epsilon})$ is big. We may write
$$m_0P^A_{\epsilon}-(K_X+\Theta^A_{\epsilon})\sim_{\mathbb Q}\sum_{i=1}^na_iE_i.$$
for some integer $n\geq 0$, irreducible curves $E_1,\dots,E_n$ on $X$, and positive real numbers $a_1,\dots,a_n$.

For any integer $1\leq i\leq n$, if $E_i\subset\Supp N$, then $(m_0P^A_{\epsilon}-(K_X+\Theta^A_{\epsilon}))\cdot E_i\geq 0$. If $E_i\not\subset\Supp N$, by Lemma \ref{lem: structure of nullA}, $E_i\not\subset\Null(P^A_{\epsilon})$, which implies that $P^A_{\epsilon}\cdot E_i>0$. Thus there exists an integer $m_1>m_0$, such that 
\begin{itemize}
\item $m_1P^A_{\epsilon}-(K_X+\Theta^A_{\epsilon})$ is nef and big, and
\item $m_1P^A_{\epsilon}$ is integral.
\end{itemize}
Since $\Supp\Theta^A_{\epsilon}\subset\Supp N$, $\Theta^A_{\epsilon}$ is simple normal crossing. By the Kawamata-Viehweg vanishing theorem, for any integer $m>0$ such that $m_1|m$, we have
$$H^i(X, mP^A_{\epsilon})=H^i(X,K_X+(mP^A_{\epsilon}-(K_X+\Theta^A_{\epsilon}))+\Theta^A_{\epsilon})=0.$$
\end{proof}

Now we state a more general vanishing theorem for foliated surface pairs:

\begin{lem}\label{lem: lem: kv vanishing foliation version 2} Assume that
\begin{itemize}
    \item $(X,\Ff,\Delta)\in\mathcal{D}_0$ is a foliated triple,
    \item $P$ and $N$ are the positive part and the negative part of the Zariski decomposition of $K_{\Ff}+\Delta$ respectively, and
    \item any irreducible component of a $(K_{\Ff}+\Delta)$-chain is not a $(-1)$-curve.
\end{itemize}
Then for any sufficiently divisible integer $m>0$ and any integer $i>0$,
$$H^i(X, mP)=0.$$
\end{lem}

\begin{proof}

\noindent\textbf{Step 1}. In this step we construct auxiliary $\mathbb Q$-divisors $\Theta$, $\Psi$ and $\Gamma$ on $X$.

 By Lemma \ref{lem: find pair for kv vanishing}, there exists a $\mathbb Q$-divisor $\Theta$ on $X$, such that
\begin{itemize}
\item $\Supp\Theta\subset\Supp N$,
\item $\lfloor\Theta\rfloor=0$,
\item for any maximal $(K_\Ff+\Delta)$-chain $C=\cup_{i=1}^{n}C_i$, 
\begin{itemize}
\item for any $1\leq i\leq n-1$, $(K_X+\Theta)\cdot C_i\leq 0$, and
\item $(K_X+\Theta)\cdot C_n\leq -1$, 
\end{itemize}
and
\item for any irreducible curve $E$ such that $E^2=-2$ and $E$ is a maximal $(K_\Ff+\Delta)$-chain, $\mu_E\Theta=\frac{1}{2}$.
\end{itemize}

By Lemma \ref{lem: structure of null}, we may define
$$\Psi:=\cup_{k=1}^r\Psi_k$$
to be the union of all the curves in $\Null(P)$ of \textbf{Type B} and \textbf{Type C}, and
$$\Gamma:=\cup_{j=1}^{s}\Gamma_j$$
to be the union of all the curves in $\Null(P)$ of \textbf{Type D}, \textbf{Type E} and \textbf{Type F}, where for any $1\leq k\leq r$ and any $1\leq j\leq s$, $\Psi_k$ and $\Gamma_j$ are connected components of $\Null(P)\backslash\Supp N$. By Lemma \ref{lem: structure of null}, each $\Psi_k$ is either a nodal curve or a cycle of smooth rational curves, and each $\Gamma_j$ is a string.\vspace{2mm}

\noindent\textbf{Step 2}. In this step we show that $h^1(\Psi\cup\Gamma,mP)=0$ for sufficiently divisible integer $m>0$.

First we show that $h^1(\Psi,mP)=0$. For any $1\leq k\leq r$, $\Psi_k$ is Cohen-Macaulay and has trivial dualizing sheaf. Moreover, by Lemma \ref{lem: structure of null}  and \cite[Theorem IV.2.2]{McQ08}, for any $1\leq k\leq r$, $K_\Ff|_{\Psi_k}$ is not torsion. By Lemma \ref{lem: Zariski decomposition and Null}(2),
$$\Psi_k\subset\Null(P)=\textbf{B}_{+}(K_\Ff+\Delta).$$
Since $\Delta\wedge\textbf{B}_{+}(K_\Ff+\Delta)=0$,
$$\Psi_k\wedge\Delta=0$$
for any $1\leq k\leq r$. By Lemma \ref{lem: structure of null}, $\Delta\cdot\Psi_k=0$ for any $1\leq k\leq r$. Thus $\Delta$ does not intersect $\Psi_k$, which implies that
$$(K_\Ff+\Delta)|_{\Psi_k}=K_\Ff|_{\Psi_k}$$ is not torsion. 

By Lemma \ref{lem: structure of null}, $\Psi_k$ does not intersect $\Supp N$. Thus for any integer $m>0$ such that $mP$ and $m(K_\Ff+\Delta)$ are integral, we have
$$\mathscr{O}_{\Psi_k}(mP)=\mathscr{O}_{\Psi_k}(m(K_\Ff+\Delta))=\mathscr{O}_{\Psi_k}(mK_\Ff).$$
By Serre duality, for any integer $m>0$ such that $mP$ and $m(K_\Ff+\Delta)$ are integral and any $1\leq k\leq r$, we have
$$h^1(\Psi_k,mP)=h^0(\Psi_k,-mP)=0,$$
which implies that
$$h^1(\Psi,mP)=0.$$

Next we show that $h^1(\Gamma,mP)=0$. For any integer $m>0$ such that $mP$ is integral, any irreducible component $C$ of $\Gamma$, we have
$$\deg(mP)|_{C}=0.$$
Since $C$ is a smooth rational curve, $g(C)=0$ where $g(C)$ is the genus of $C$. Thus
$$0=\deg(mP)|_{C}\geq 2g(C)-2.$$
By Serre duality, 
$$h^1(C, mP)=0,$$
which implies that
$$h^1(\Gamma,mP)=0.$$
Thus $h^1(\Psi\cup\Gamma,mP)=0.$\vspace{2mm}

\noindent\textbf{Step 3}. In this step we use exact sequences to reduce the vanishing of $h^i(X, mP)$ to the vanishing of $h^i(X,mP-\Psi-\Gamma)$.

Since $\Psi\cap\Gamma=\emptyset$, for any integer $m>0$ such that $mP$ is integral, we have the short exact sequence
$$0\rightarrow \mathscr{O}_{X}(mP-\Psi-\Gamma)\rightarrow \mathscr{O}_{X}(mP)\rightarrow \mathscr{O}_{\Psi\cup\Gamma}(mP)\rightarrow 0,$$
which induces the exact sequence
\begin{align*}
0\rightarrow & H^0(X, mP-\Psi-\Gamma)\rightarrow H^0(X, mP)\rightarrow H^0(\Psi\cup\Gamma, mP)\\
\rightarrow & H^1(X, mP-\Psi-\Gamma)\rightarrow H^1(X, mP)\rightarrow 0
\end{align*}
and
$$H^2(X, mP-\Psi-\Gamma)\cong H^2(X, mP).$$
Therefore, in order to show that $h^i(X, mP)=0$ for any integer $i>0$ and sufficiently divisible integer $m>0$, we only need to show that $h^i(X, mP-\Psi-\Gamma)=0$ for any integer $i>0$ and sufficiently divisible integer $m>0$.\vspace{2mm}

\noindent\textbf{Step 4}. Let $Q:=K_X+\Theta+\Psi+\Gamma$. For each irreducible curve $C\subset\Null(P)$, we show that $Q\cdot C\leq 0$.  By Lemma \ref{lem: structure of null}, there are six cases.\vspace{2mm}

\noindent\textbf{Case 1}. $C$ is of \textbf{Type A}. In this case, $\Psi\cdot C=0$. There are two possibilities:\vspace{2mm}

\noindent\textbf{Case 1.1} $C\cap \Gamma=\emptyset$. By Lemma \ref{lem: find pair for kv vanishing}(3), 
$$Q\cdot C=(K_X+\Theta)\cdot C\leq 0.$$ 

\noindent\textbf{Case 1.2} $C\cap\Gamma\not=\emptyset$. By Lemma \ref{lem: structure of zdecomposition of foliation}(1), $C$ is a part of a $(K_\Ff+\Delta)$-chain. Since $\Gamma$ does not contain any component of any $(K_\Ff+\Delta)$-chain and each irreducible component of $\Gamma$ is $\Ff$-invariant, $C$ is the last curve of a maximal $(K_\Ff+\Delta)$-chain $D$ and $\Gamma$ contains the tail of $D$. Therefore $\Gamma\cdot C=1$, and by Lemma \ref{lem: find pair for kv vanishing}(3), $(K_X+\Theta)\cdot C\leq -1$. Thus
$$Q\cdot C=(K_X+\Theta)\cdot C+\Gamma\cdot C\leq -1+1=0.$$

\noindent\textbf{Case 2}. $C$ is of \textbf{Type B}. In this case, $C\cap\Supp N=\emptyset$, $\Gamma\cap C=\emptyset$, and $\Psi\cdot C=C^2$.  Thus 
$$Q\cdot C=(K_X+\Psi)\cdot C=(K_X+C)\cdot C=-\chi(C)=0.$$

\noindent\textbf{Case 3}. $C$ is of \textbf{Type C}.  In this case, $C\cap\Supp N=\emptyset$, $\Gamma\cap C=\emptyset$, and $\Psi\cdot C=C^2+2$. Thus 
$$Q\cdot C=(K_X+\Psi)\cdot C=(K_X+C)\cdot C+2=-\chi(C)+2=0.$$

\noindent\textbf{Case 4}. $C$ is of \textbf{Type D}. In this case, $C\cap\Supp N=\emptyset$ and $\Psi\cap C=\emptyset$. Since the connected component of $\Gamma$ containing $C$ is a string, $\Gamma\cdot C\leq C^2+2.$ Thus 
$$Q\cdot C=(K_X+\Gamma)\cdot C\leq (K_X+C)\cdot C+2=-\chi(C)+2=0.$$

\noindent\textbf{Case 5}. $C$ is of \textbf{Type E}. In this case, there is a $(K_{\Ff}+\Delta)$-chain $D=\cup_{i=1}^nC_i$ such that $D\cup C$ is an $\Ff$-chain, and $\Psi\cap C=\emptyset$. Since the coefficients of $\Theta$ belong to $[0,1)$, $\Supp\Theta\subset\Supp N$, and $C_n$ is the only irreducible curve of $\Supp N$ which intersects $C$, we have $\Theta\cdot C<1.$ Since $\Gamma$ only contains $\Ff$-invariant curves and $D\cup C$ is an $\Ff$-chain, $C$ intersects at most one other irreducible curve which belongs to $\Gamma$, which implies that $\Gamma\cdot C\leq C^2+1.$ We deduce that
\begin{align*}
Q\cdot C&=(K_X+\Theta+\Gamma)\cdot C<(K_X+\Gamma)\cdot C+1\\
&\leq (K_X+C)\cdot C+2=-\chi(C)+2=0.
\end{align*}

\noindent\textbf{Case 6}. $C$ is of \textbf{Type F}. Then $\Psi\cap C=\emptyset$ and $C$ intersects at most one other irreducible curve which belongs to $\Gamma$. Thus $\Gamma\cdot C\leq C^2+1.$ Moreover, by our assumptions, for any $(K_{\Ff}+\Delta)$-chain $D$ such that $\det(-||D||)=2$, $D$ is an irreducible curve with self-intersection $-2$. By the construction of $\Theta$, $\Theta\cdot C=1.$ We deduce that
$$Q\cdot C=(K_X+\Theta+\Gamma)\cdot C\leq (K_X+C)\cdot C+2=-\chi(C)+2=0.$$

\noindent\textbf{Step 5}. Since $K_{\Ff}+\Delta$ is big, $P$ is big. Thus there exists a sufficiently divisible integer $m_0>0$ such that
\begin{itemize}
\item $m_0P-Q$ is big, and
\item $m_0P$ is integral.
\end{itemize}
In particular, we may write
$$m_0P-Q\sim_{\mathbb Q}\sum_{i=1}^la_iE_i$$
for some integer $l\geq 0$, irreducible curves $E_1,\dots,E_l$ on $X$, and positive real numbers $a_1,\dots,a_l$.

For any integer $1\leq i\leq l$, if $P\cdot E_i=0$, then $Q\cdot C\leq 0$, which implies that $(m_0P-Q)\cdot C\geq 0$. Thus there exists an integer $m_1>0$, such that $m_0|m_1$, and $(m_1P-Q)\cdot E_i\geq 0$ for any $E_i$. In particular, $m_1P-Q$ is big and nef, hence $mP-Q$ is big and nef for any integer $m>0$ such that $m_1|m$.

Since $\Supp\Theta\subset\Supp N$, $\Theta$ is simple normal crossing. By the Kawamata-Viehweg vanishing theorem, for any integer $m>0$ such that $m_1|m$ and $i>0$, 
$$H^i(X, mP-\Psi-\Gamma)=H^i(X,K_X+(mP-(K_X+Q)+\Theta))=0,$$
and the lemma follows.
\end{proof}

\begin{proof}[Proof of Theorem \ref{thm: weak kv vanishing foliated surface}] We first show that $\Delta\wedge\textbf{B}_{+}(K_\Ff+\Delta)=0$. Suppose not, then there exists an irreducible curve $C\subset\Supp\Delta$ such that $C\subset\textbf{B}_{+}(K_\Ff+\Delta)$. Since $K_\Ff+\Delta$ is big and nef and $(K_{\Ff}+\Delta)\cdot C=0$, we have $C^2<0$. Since $\lfloor\Delta\rfloor=0$, we have
$$\Delta\cdot C\geq(\mult_C\Delta)\cdot C^2>C^2.$$ 
Since $C\subset\Supp\Delta$, $C$ is not $\Ff$-invariant. By Lemma \ref{lem: tang formula},
$$0\leq \tang(\Ff,C)=(K_{\Ff}+C)\cdot C<(K_{\Ff}+\Delta)\cdot C=0,$$
a contradiction. Thus $\Delta\wedge\textbf{B}_{+}(K_\Ff+\Delta)=0$. 

Since $K_\Ff+\Delta$ is nef, $K_{\Ff}+\Delta$ is the positive part of the Zariski decomposition of itself, and Theorem \ref{thm: weak kv vanishing foliated surface} follows from Lemma \ref{lem: lem: kv vanishing foliation version 2}.\end{proof}

\section{Deformation theory of foliated surface pairs}

In this section we prove several results on deformation of foliated surface pairs. The main result of this section shows that the negative part of the Zariski decomposition of certain foliated surface pairs deforms, which is one of the key ingredient to prove Theorem \ref{thm: foliated surface pair inv plur}.

We first state a lemma on intersection numbers that will be repeatedly used in the rest of this section:
\begin{lem}\label{lem: inter deform imply kfb chain deform}
Let $n>0$ be an integer. For any $i\in\{1,2\}$, assume that
\begin{itemize}
    \item $X_i\in\mathcal{S}_{cyc}$ is a surface,
    \item $\Ff_i$ a canonical foliation on $X_i$, 
    \item $\Delta_i$ is an effective $\mathbb Q$-divisor on $X_i$, 
    \item  $C_i:=\cup_{j=1}^nC_{j,i}$ is an $\Ff_i$-chain, such that
    \begin{itemize}
        \item $C_1$ is a $(K_{\Ff_1}+\Delta_1)$-chain,
        \item for any $1\leq j\leq n$, $(K_{\Ff_1}+\Delta_1)\cdot C_{j,1}=(K_{\Ff_2}+\Delta_2)\cdot C_{j,2}$, and
        \item for any $1\leq j\leq n$, $C_{j,1}^2=C_{j,2}^2$.
    \end{itemize}
\end{itemize}
Then $C_2$ is a $(K_{\Ff_2}+\Delta_2)$-chain. Moreover, if $X_1,X_2$ are smooth, and $C_1$ is a $(K_{\Ff_1}+\Delta_1)$-artificial chain, then $C_2$ is a $(K_{\Ff_2}+\Delta_2)$-artificial chain.
\end{lem}
\begin{proof}
Since
$$(K_{\Ff_2}+\Delta_2)\cdot C_{1,2}=(K_{\Ff_1}+\Delta_1)\cdot C_{1,1}<0,$$
$C_{1,2}$ is a $(K_{\Ff_2}+\Delta_2)$-negative extremal ray. Thus there exists a contraction $\phi_2: (X_2,\Ff_2,\Delta_2)\rightarrow (X_2',\Ff_2',\Delta'_2)$ of $C_{1,2}$ and a contraction $\phi_1:(X_1,\Ff_1,\Delta_1)\rightarrow (X_1',\Ff_1',\Delta'_1)$ of $C_{1,1}$. 

If $n=1$, we are done. Otherwise, for any $i\in\{1,2\}$, let $C_i'=(\phi_i)_*C_i$ and $C_{j,i}':=(\phi_i)_*C_{j,i}$ for any $2\leq j\leq n$. Then 
\begin{itemize}
    \item $C_i'$ is a $\Ff'_i$-chain for any $i\in\{1,2\}$,
    \item $C_1'$ is a $(K_{\Ff_1'}+\Delta_1')$-chain,
    \item for any $2\leq j\leq n$, $(K_{\Ff_1'}+\Delta_1')\cdot C'_{j,1}=(K_{\Ff'_2}+\Delta'_2)\cdot C'_{j,2}$, and
\item for any $2\leq j\leq n$, $(C'_{j,1})^2=(C'_{j,2})^2$.
\end{itemize} 
Possibly replacing $X_i,\Ff_i,\Delta_i$ and $C_i$ with $X_i',\Ff_i',\Delta_i'$ and $C_i'$ respectively and by applying induction on $n$, we deduce that $C_2$ is a $(K_{\Ff_2}+\Delta_2)$-chain. 

If $X_1,X_2$ are smooth and $C_1$ is a $(K_{\Ff_1}+\Delta_1)$-artificial chain, by Lemma \ref{lem: contract artificial chain step by step}, there exists a sequence of contractions of $(-1)$-curves
$$(X_{1,1},\Ff_{1,1},\Delta_{1,1}):=(X_1,\Ff_1,\Delta_1)\xrightarrow{\phi_1}\dots\xrightarrow{\phi_n}(X_{n,1},\Ff_{n,1},\Delta_{n,1}),$$
such that $\phi_n\circ\dots\phi_1$ is the contraction of $C_1$. For every $1\leq i\leq n$, suppose that $\phi_i$ contracts the strict transform of $C_{r_i,1}$ on $X_{i-1,1}$ for some $1\leq r_i\leq n$. We inductively define a sequence of contractions
$$(X_{1,2},\Ff_{1,2},\Delta_{1,2}):=(X_2,\Ff_2,\Delta_2)\xrightarrow{\phi'_1}\dots\xrightarrow{\phi'_n}(X_{n,2},\Ff_{n,2},\Delta_{n,2}),$$
such that each $\phi_i'$ is the contraction of the strict transform of $C_{r_i,2}$ on $X_{i-1,2}$. Then each $\phi_i'$ is a contraction of an $\Ff_{i-1,2}$-invariant curve with self-intersection $-1$. The lemma follows from Lemma \ref{lem: contraction changes z}.
\end{proof}

\begin{lem}\label{lem: kf chain deforms}  Assume that
\begin{itemize}
    \item $X$ is a smooth threefold,
    \item $\Ff$ is a rank 1 foliation on $X$,
    \item $\pi: X\rightarrow T$ is a smooth morphism to a curve, such that $(X_t,\Ff_t)|_{t\in T}$ is a smooth family of canonical foliations of surfaces associated with $\pi$.
    \item $s\in T$ is a closed point, 
    \item $n>0$ an integer, and
    \item $C_s=\cup_{i=1}^{n}C_{i,s}$ is a $K_{\Ff_s}$-chain on $X_s$,
\end{itemize}
then there exists an open neighborhood $U$ of $s$ in $T$, hypersurfaces $C_1,\dots, C_n\subset \pi^{-1}(U)$, such that  for every $1\leq i\leq n$,
\begin{enumerate}
\item $C_i|_{X_s}=C_{i,s}$,
\item for any $t\in U$,
$$C_t:=\cup_{i=1}^{n}(C_i\cap X_t)$$
is a $K_{\Ff_t}$-chain.
\end{enumerate}
\end{lem}

\begin{proof}

We use induction on $n$. When $n=1$, if $C_{i,s}^2=-1$, the lemma follows from Lemma \ref{lem: -1 curve deforms}. If $C_{i,s}^2\leq -2$, the lemma follows from Lemma \ref{lem: F chain deforms} and Lemma \ref{lem: inter deform imply kfb chain deform}.

Suppose we are done with the cases when $C_s$ has length $1,2,\dots n-1$. When $C_s$ has length $n$, let $k\geq 0$ be an integer and $1\leq r_1<\dots<r_k\leq n$ all the integers such that $(C_{r_i,s})^2=-1$ for each $i$. If $k=0$, the lemma follows from Lemma \ref{lem: F chain deforms}, so we may suppose that $k>0$. By Lemma \ref{lem: no continue -1 curve}, for any $1\leq i,j\leq k$ such that $i\not=j$, $C_{r_i,s}\cap C_{r_j,s}=\emptyset$. By Lemma \ref{lem: -1 curve deforms}, there exists an open neighborhood $U$ of $s$ in $T$, hypersurfaces $\{C_{r_i}\}_{i=1}^{k}$ over $U$, such that
\begin{itemize}
\item for any $1\leq i\leq k$, $C_{r_i}\cap X_s=C_{r_i,s}$, 
\item for any $t\in U$, $C_{r_i,t}:=C_{r_i}\cap X_t$ is a smooth rational curve with self-intersection $-1$, and
\item there is a birational contraction $\nu: \pi^{-1}(U)\rightarrow X'$ and a smooth morphism $\pi': X'\rightarrow U$, such that for every closed point $t\in U$, $\nu_t:=\nu|_{X_t}$ is the contraction of $C_{r_i,t}$.
\end{itemize}
For any $l\not\in\{r_1,\dots,r_k\}$ and any $t\in U$, let $C_{l,t}'$ be the birational transform of $C_{l,t}$ on $X'$. By Lemma \ref{lem: contraction changes z} and Lemma \ref{lem: Intersection surface}, $\cup_{l\not\in\{r_1,\dots,r_k\}}C_{l,s}'$ is a $K_{(\nu_s)_*\Ff_s}$-chain. By induction on $n$, possibly shrinking $U$, we may assume that $\cup_{l\not\in\{r_1,\dots,r_k\}}C_{l,t}'$ is a $K_{(\nu_t)_*\Ff_t}$-chain for every $t\in U$. Since $\nu_t$ is a sequence of blow-downs, the lemma follows from Lemma \ref{lem: blow up changes z} and Lemma \ref{lem: Intersection surface}.
\end{proof}

\begin{lem}\label{lem: foliated pair chain deforms} 
Assume that 
\begin{itemize}
    \item $X$ is a smooth threefold, 
    \item $\Ff$ is a rank 1 foliation on $X$,
    \item $\pi: X\rightarrow T$ is a smooth morphism to a curve,
    \item $s\in T$ is a closed point, and
    \item $\Delta$ is an effective $\mathbb Q$-Cartier $\mathbb Q$-divisor on $X$
\end{itemize}
satisfying the following:
\begin{itemize}
    \item $(X_t,\Ff_t)|_{t\in T}$ is a smooth family of canonical foliations of surfaces associated with $\pi$,
    \item $\Delta_t:=\Delta|_{X_t}$ for any closed point $t\in T$,
    \item for any closed point $t\in T$ and any irreducible component $\Delta^0$ of $\Delta$, $\Delta^0|_{X_t}$ is irreducible and reduced, and 
    \item $(X_s,\Ff_s,\Delta_s)\in\mathcal{D}$.
\end{itemize}
Suppose that $C_s=\cup_{i=1}^{n}C_{i,s}$ is a $(K_{\Ff_s}+\Delta_s)$-chain. Then there exists an open neighborhood $U$ of $s$ in $T$ and hypersurfaces $E_1,\dots,E_n\subset\pi^{-1}(U)$ which are transversal to the fibers of $\pi$, such that for any $1\leq i\leq n$ and any $t\in U$,
\begin{enumerate}
\item $E_i\cap X_s=C_{i,s}$,
\item $C_{i,t}:=E_i\cap X_t$ is a smooth rational curve,
\item $C_t:=\cup_{i=1}^{n}(E_{i}\cap X_t)$ is a $(K_{\Ff_t}+\Delta_t)$-chain, and
\item if $C_s$ is a $(K_{\Ff_s}+\Delta_s)$-artificial chain, then $C_t$ is a $(K_{\Ff_t}+\Delta_t)$-artificial chain.
\end{enumerate}
\end{lem}

\begin{proof} By Lemma \ref{lem: kfb chain is kf chain}, $C_s$ is a $K_{\Ff_s}$-chain. By Lemma \ref{lem: kf chain deforms}, there exists an open neighborhood $U$ of $s$ and hypersurfaces $E_1,\dots,E_n\subset \pi^{-1}(U)$ which are transversal to the fibers of $\pi$, such that \begin{itemize}
    \item intersection numbers deform over $U$,
    \item for any $t\in U$, $C_{i,t}:=E_i\cap X_t$ is a smooth rational curve, and
    \item $C_t:=\cup_{i=1}^{n}C_{i,t}$ is a $K_{\Ff_t}$-chain. 
\end{itemize}
The lemma follows from Lemma \ref{lem: inter deform imply kfb chain deform}.
\end{proof}

\begin{lem}\label{lem: inv plur foliated pairs under bir contraction}  
Assume that 
\begin{itemize}
    \item $X$ is a smooth threefold, 
    \item $\Ff$ is a rank 1 foliation on $X$,
    \item $\pi: X\rightarrow T$ is a smooth morphism to a curve,
    \item $s\in T$ is a closed point,
    \item $\Delta$ is an effective $\mathbb Q$-Cartier $\mathbb Q$-divisor on $X$, and
    \item $m>0$ is an integer
\end{itemize}
satisfying the following:
\begin{itemize}
    \item  $(X_t,\Ff_t)|_{t\in T}$ is a smooth family of canonical foliations of surfaces associated with $\pi$,
    \item $m\Delta$ is integral,
    \item $\Delta_t:=\Delta|_{X_t}$ for any closed point $t\in T$,
    \item for any closed point $t\in T$ and any irreducible component $\Delta^0$ of $\Delta$, $\Delta^0|_{X_t}$ is irreducible and reduced, 
    \item $K_{\Ff}+\Delta$ is $\pi$-pseudo-effective,
    \item $\Delta_s\wedge \textbf{\rm \textbf{B}}_{-}(K_{\Ff_s}+\Delta_s)=0$, and
    \item there is an irreducible curve $D_s\subset X_s$, such that
    \begin{itemize}
        \item $D_s^2=-1$, and
        \item $D_s$ is an irreducible component of a $(K_{\Ff_s}+\Delta_s)$-chain.
    \end{itemize}
\end{itemize}
Then there is
\begin{itemize}
\item an open subset $s\in U$ of $T$,
\item a birational contraction $f: (\pi^{-1}(U),\Ff|_{\pi^{-1}(U)},\Delta|_{\pi^{-1}(U)})\rightarrow (X',\Ff',\Delta')$ over $U$ and
\item a smooth morphism $\pi':X'\rightarrow U$, 
\end{itemize}
satisfying the following.  For any $t\in U$, let $f_t:=f|_{\pi^{-1}(t)}$, $X'_t:=(\pi')^{-1}(t)$, $\Ff'_t:=\Ff'|_{X'_t}$ and $\Delta_t':=\Delta'|_{X'_t}$, then
\begin{enumerate}
\item $(X_t',\Ff_t')|_{t\in U}$ is a smooth family of canonical foliations of surfaces, 
\item $f_t$ is a birational contraction for any $t\in U$,
\item $f_s$ is the contraction of $D_s$, 
\item $\pi'\circ f=\pi|_{\pi^{-1}(U)}$,
\item $K_{\Ff'}+\Delta'$ is $\pi'$-pseudo-effective, 
\item $\Delta'_s\wedge \textbf{\rm \textbf{B}}_{-}(K_{\Ff'_s}+\Delta'_s)=0$, and
\item $$H^0(X_t',m(K_{\Ff_t'}+\Delta'_t))\cong H^0(X_t,m(K_{\Ff_t}+\Delta_t))$$ for any $t\in U$.
\end{enumerate}
\end{lem}

\begin{proof} Since $D_s$ is an irreducible component of a $(K_{\Ff_s}+\Delta_s)$-chain, there exists a $(K_{\Ff_s}+\Delta_s)$-chain $C_s=\cup_{i=1}^{n}C_{i,s}$ such that $D_s=C_{n,s}$.  By Lemma \ref{lem: -1 curve deforms} and Lemma \ref{lem: foliated pair chain deforms}, there exists an an open neighborhood $U$ of $s$, hypersurfaces $E_1,\dots,E_n\subset\pi^{-1}(U)$ which are transversal to the fibers of $\pi$, such that 
\begin{itemize}
\item for any $1\leq i\leq n$, $E_i\cap X_s=C_{i,s}$, and 
\item for any $t\in U$, 
\begin{itemize}
\item $C_{i,t}:=E_i|_{X_t}$ is a smooth rational curve, 
\item $C_t:=\cup_{i=1}^{n}C_{i,t}$ is a $(K_{\Ff_t}+\Delta_t)$-chain, and
\item $C_{n,t}^2=-1$.
\end{itemize}
\end{itemize}
Let
$$f: \pi^{-1}(U)\rightarrow X'$$
be the contraction of $E_n$ and $\pi':X'\rightarrow U$ the induced morphism, then (1)-(6) follow immediately, and we only need to show (7). 

For any $t\in U$, we define $g_t$ in the following way:
\begin{itemize}
\item if $n=1$, let $g_t:=\id_{X'_t}$, and 
\item if $n>1$, let $C'_t:=(f_t)_*C_t$ for any $1\leq i\leq n-1$. By applying  Lemma \ref{lem: Intersection surface}, $C'_t$ is a $(K_{\Ff_t'}+\Delta'_t)$-chain. We let $g_t$ be the contraction of $C'_t$.
\end{itemize}
Then $g_t\circ f_t$ is the contraction of $C_t$. Thus
\begin{align*}
H^0(X_t',m(K_{\Ff_t'}+\Delta'_t))&\cong H^0(X_t',m(K_{(g_t)_*{\Ff'_t}}+(g_t)_*\Delta'_t))\\
&=H^0(X_t',m(K_{(g_t\circ f_t)_*{\Ff_t}}+(g_t\circ f_t)_*\Delta_t))\\
&\cong H^0(X_t,m(K_{\Ff_t}+\Delta_t)),
\end{align*}
which implies (7).\end{proof}

\begin{lem}\label{lem: deform zariski decomposition foliated pairs} 
Assume that 
\begin{itemize}
    \item $X$ is a smooth threefold, 
    \item $\Ff$ is a rank 1 foliation on $X$,
    \item $\pi: X\rightarrow T$ is a smooth morphism to a curve,
    \item $s\in T$ is a closed point, and
    \item $\Delta$ is an effective $\mathbb Q$-Cartier $\mathbb Q$-divisor on $X$
\end{itemize}
satisfying the following:
\begin{itemize}
    \item  $(X_t,\Ff_t)|_{t\in T}$ is a smooth family of canonical foliations of surfaces associated with $\pi$,
    \item $\Delta_t:=\Delta|_{X_t}$ for any closed point $t\in T$,
    \item for any closed point $t\in T$ and any irreducible component $\Delta^0$ of $\Delta$, $\Delta^0|_{X_t}$ is irreducible and reduced, 
\item $K_{\Ff}+\Delta$ is $\pi$-pseudo-effective, 
\item $\Delta_s\wedge\textbf{\rm \textbf{B}}_{-}(K_{\Ff_s}+\Delta_s)=0$, and
\item $X_s$ does not contain any $(K_{\Ff_s}+\Delta_s)$-artificial chain.
\end{itemize}
Then there is an open neighborhood $U$ of $s$ in $T$ and an effective $\Qq$-divisor $N\subset\pi^{-1}(U)$, such that for every $t\in U$,
\begin{enumerate}
    \item $N|_{X_t}$ is the negative part of the Zariski decomposition of $K_{\Ff_t}+\Delta_t$, and
    \item for every irreducible component $D$ of $N$, $D$ intersects $X_t$ transversally and $D\cap X_t$ is a smooth rational curve.
\end{enumerate}
\end{lem}

\begin{proof} The proof follows from the same lines as in \cite[Proposition 1]{Bru01} and \cite[Proposition 2]{Bru01}. For readers' convenience, we give a full proof here. \vspace{2mm}

For every $t\in T$, suppose that $P_t$ and $N_t$ are the positive part and the negative part of the Zariski decomposition of $K_{\Ff_t}+\Delta_t$ respectively. By Lemma \ref{lem: structure of zdecomposition of foliation}, $\Supp N_s$ is the union of all the maximal $(K_{\Ff_s}+\Delta_s)$-chains. By Lemma \ref{lem: foliated pair chain deforms}, there is an open neighborhood $U_0$ of $s$ and an effective $\mathbb Q$-divisor $M\subset\pi^{-1}(U_0)$, such that 
\begin{itemize}
\item $M|_{X_s}=N_s$, 
\item for every $t\in U_0$ and every irreducible component $D$ of $M$, $D$ meets $X_t$ transversally and $D\cap X_t$ is a smooth rational curve, and
\item for every $t\in U_0$, every $(K_{\Ff_s}+\Delta_s)$-chain which belongs to $\Supp N_s$ deforms to a $(K_{\Ff_t}+\Delta_t)$-chain which belongs to $\Supp M|_{X_t}$, 
\end{itemize}
We let $M_t:=M|_{X_t}$ for any $t\in U_0$. In the rest of the proof, we show that there is an open neighborhood $s\in U\subset U_0$, such that $N_t=M_t$ for any $t\in U$.\vspace{2mm}

\begin{claim}\label{claim: suppmsubsetsuppn}
$\Supp M_t\subset\Supp N_t$ for any $t\in U_0$.
\end{claim}

\begin{proof}[Proof of Claim \ref{claim: suppmsubsetsuppn}]
Let $C_s=\cup_{i=1}^nC_{i,s}$ be a maximal $(K_{\Ff_s}+\Delta_s)$-chain. For any $t\in U_0$, suppose that $C_s$ deforms to $C_t$ and $C_{i,s}$ deforms to $C_{i,t}$ for any $1\leq i\leq n$. Then each $C_t$ is a $(K_{\Ff_t}+\Delta_t)$-chain. 
In particular, there is a sequence of birational contractions
$$(X_t^0,\Ff_t^0,\Delta_t^0):=(X_t,\Ff_t,\Delta_t)\xrightarrow{\phi_t^1}(X_t^1,\Ff_t^1,\Delta_t^1)\dots\xrightarrow{\phi_t^n}(X_t^n,\Ff_t^n,\Delta_t^n),$$
such that for every $1\leq i\leq n$,
\begin{itemize}
    \item $(X_t^i,\Ff_t^i,\Delta_t^i)\in\mathcal{D}$,
    \item $\phi_{t}^i$ is the contraction of a $(K_{\Ff_t^{i-1}}+\Delta_t^{i-1})$-negative extremal ray $C_{i,t}'$ on $X_t^{i-1}$, and
    \item $C_{i,t}'$ is the strict transform of $C_{i,t}$ on $X_t^{i-1}$.
\end{itemize}
We use induction on $i$ to show that $C_t\subset\Supp N_t$. For any integer $1\leq i\leq n$, suppose that $C_{1,t},\dots,C_{i-1,t}\subset\Supp N_t$. Let $\nu:=\id_{X_t}$ when $i=1$ and $\nu:=\phi_{t}^{i-1}\circ\dots\phi_t^1$ when $i>1$. Then
$$0>(K_{\Ff_t^{i-1}}+\Delta_t^{i-1})\cdot C'_{i,t}=\nu_*(K_{\Ff_t}+\Delta_t)\cdot C'_{i,t}=\nu_*(P_t+N_t)\cdot C'_{i,t}.$$
Since $P_t$ is nef, by Lemma \ref{lem: Intersection surface}, $\nu_*P_t$ is nef. Thus $\nu_*N_t\cdot C_{i,t}'<0$, which implies that $C_{i,t}'\subset\Supp(\nu_*N_t)$. Thus $C_{i,t}\subset\Supp N_t$. By induction on $i$, $C_t\subset\Supp N_t$.
\end{proof}

\begin{claim}\label{claim: ntgeqmt} There is an open neighborhood $U_1\subset U_0$ of $s$, such that for any $t\in U_1$,
\begin{enumerate} 
    \item $N_t\geq M_t$, and
    \item for any irreducible curve $L_t\subset\Supp M_t$, $N_t\cdot L_t=M_t\cdot L_t$.
\end{enumerate}
\end{claim}

\begin{proof}[Proof of Claim \ref{claim: ntgeqmt}]
Let $L_t\subset\Supp M_t$ be an irreducible curve. Since $\Supp M_t\subset\Supp N_t$, $P_t\cdot L_t=0$. Let $U_1\subset U_0$ be the locus such that intersection numbers deform over $U_1$. Since $L_t\subset\Supp M_t$, there is an irreducible curve $L_s\subset\Supp N_s$ such that $L_t$ is deformed by $L_s$. Thus
\begin{align*}
N_t\cdot L_t&=(P_t+N_t)\cdot L_t=(K_{\Ff_t}+\Delta_t)\cdot L_t=(K_{\Ff_s}+\Delta_s)\cdot L_s\\
&=(P_s+N_s)\cdot L_s=N_s\cdot L_s=M_t\cdot L_t,
\end{align*}
which implies (2).

Let $\nu_t$ be the contraction of $\Supp M_t$. By (2), $M_t-N_t$ is $\nu_t$-nef. Since $\Supp M_t\subset\Supp N_t$, we may write 
$$M_t-N_t=A_t-B_t,$$
where $B_t\geq 0$, $\Supp B_t\cap\Supp M_t=0$, and $\Supp A_t\subset\Supp M_t$. Then $A_t$ is $\nu_t$-nef. By Lemma \ref{lem: num negativity lemma}, $-A_t\geq 0$, which implies (1).
\end{proof}

%
%
\noindent\textit{Proof of Lemma \ref{lem: deform zariski decomposition foliated pairs} continued}. Let $Q_t:=P_t+(N_t-M_t)$ for every $t\in U_1$. 

For every $t\in U_1$ such that $Q_t$ is nef, by Claim \ref{claim: ntgeqmt}(1), $N_t\geq M_t$. Since the intersection matrix of $N_t$ is negative definite, the intersection matrix of $N_t-M_t$ is negative definite. Moreover, for any curve $L_t\subset\Supp(N_t-M_t)\subset\Supp N_t$, we have $P_t\cdot (N_t-M_t)=0$. Thus $Q_t=Q_t+0$ and $Q_t=P_t+(N_t-M_t)$ are both the Zariski decomposition of $Q_t$, which implies that $N_t=M_t$. 

For every $t\in U_1$ such that $Q_t$ is not nef, there is an irreducible curve $G_t$ on $X_t$ such that $Q_t\cdot G_t<0$. Since $P_t$ is nef and $N_t-M_t\geq 0$,  $G_t\subset\Supp N_t$. Thus $P_t\cdot G_t=0$ and $0>Q_t\cdot G_t=(N_t-M_t)\cdot G_t$. By Claim \ref{claim: ntgeqmt}(2), $G_t\not\subset\Supp M_t$. 

\begin{claim}\label{claim: qt not nef} For every $t\in U_1$ such that $Q_t$ is not nef and every irreducible curve $G_t$ on $X_t$ such that $Q_t\cdot G_t<0$,
\begin{enumerate}
    \item $G_t\not\subset\Supp\Delta_t$, and
    \item $G_t$ is $\Ff_t$-invariant.
\end{enumerate}
\end{claim}
\begin{proof}[Proof of Claim \ref{claim: qt not nef}]
Suppose that $G_t\subset\Supp\Delta_t$, then $G_t$ is deformed from an irreducible curve $G_s\subset\Supp\Delta_s\subset X_s$. Thus
$$0=P_s\cdot G_s=Q_s\cdot G_s=Q_t\cdot G_t<0,$$
a contradiction, and we deduce (1).

Since $M_t$ is a union of $(K_{\Ff_t}+\Delta_t)$-chains, we may let $\phi_t: (X_t,\Ff_t,\Delta_t)\rightarrow (X_t',\Ff_t',\Delta_t')$ be the contraction of all the connected components of $\Supp M_t$ which intersect $G_t$, $G_t':=(\phi_t)_*G_t$ and $Q_t':=(\phi_t)_*Q_t$. By Lemma \ref{lem: Intersection surface}, 
$$Q'_t\cdot G'_t=Q_t\cdot G_t<0.$$
Since $P_t$ is nef and $N_t\geq M_t$, $Q_t$ is pseudo-effective. Thus $Q_t'$ is pseudo-effective, which implies that $(G_t')^2<0.$ By Lemma \ref{lem: tang formula}, if $G_t'$ is not $\Ff_t'$-invariant, then $K_{\Ff_t'}\cdot G_t'>0$.

Since
\begin{align*}
0>Q'_t\cdot G'_t&=(\phi_t)_*(P_t+N_t-M_t)\cdot G_t'=(\phi_t)_*(P_t+N_t)\cdot G_t'\\
&=(\phi_t)_*(K_{\Ff_t}+\Delta_t)\cdot G_t'=(K_{\Ff'_t}+\Delta'_t)\cdot G_t',
\end{align*}
if $G_t'$ is not $\Ff_t'$-invariant, then $\Delta'_t\cdot G_t'<0$. This implies that $G'_t\subset\Supp\Delta'_t$, which contradicts to (1). Thus $G_t'$ is not $\Ff_t'$-invariant, which implies (2).
\end{proof}

\begin{claim}\label{claim: gt is kft chain}
For every $t\in U_1$ such that $Q_t$ is not nef and every irreducible curve $G_t$ on $X_t$ such that $Q_t\cdot G_t<0$, $G_t$ belongs to a $K_{\Ff_t}$-chain.
\end{claim}
\begin{proof}[Proof of Claim \ref{claim: gt is kft chain}]
First we suppose that $G_t\cap\Supp M_t=\emptyset$. Then $N_t\cdot G_t<0$, which implies that $(K_{\Ff_t}+\Delta_t)\cdot G_t<0$. Since $K_{\Ff_t}+\Delta_t$ is pseudo-effective, $G_t^2<0$. Since $G_t\not\subset\Supp\Delta_t$, $\Delta_t\cdot G_t>0$. Thus $K_{\Ff_t}\cdot G_t<0$, which implies that $G_t$ is a $K_{\Ff_t}$-chain.

Suppose that $G_t\cap\Supp M_t\not=\emptyset$. Then there is a connected component $H_t$ of $\Supp M_t$, such that $H_t=\cup_{i=1}^lH_{i,t}$ is a $(K_{\Ff_t}+\Delta_t)$-chain and $G_t$ is the tail of $H_t$. Let $\phi_t: (X_t,\Ff_t,\Delta_t)\rightarrow (X_t',\Ff_t',\Delta_t')$ be the contraction of all the connected components of $\Supp M_t$ which intersect $G_t$, $G'_t:=(\phi_t)_*G_t$ and $Q_t':=(\phi_t)_*Q_t$. Then
$$0>Q_t\cdot G_t=Q'_t\cdot G'_t=(\phi_t)_*(K_{\Ff_t}+\Delta_t-M_t)\cdot G_t'=(K_{\Ff'_t}+\Delta'_t)\cdot G_t'.$$
Since $G_t\not\subset\Supp\Delta_t$, $G_t'\not\subset\Supp\Delta_t'$. Thus $K_{\Ff_t'}\cdot G'_t<0$. By Lemma \ref{lem: z and cs formula}(1), $\chi(G_t')>Z(\Ff_t',G_t')$. Thus $G_t'$ is a smooth rational curve. Since $K_{\Ff'_t}+\Delta'_t$ is pseudo-effective, $(G_t')^2<0.$ By Lemma \ref{lem: z and cs formula}(2), $CS(\Ff_t',G_t')<0$. Thus $G_t'\cap\Sing(\Ff_t')\not=\emptyset$, hence $Z(\Ff_t',G_t')\geq1$. Suppose that $\{p_1,\dots,p_m\}=G_t'\cap\Sing(X_t')$ such that each $p_i$ is of order $r_i\geq 2$, then
$\chi(G_t')=2-\sum_{i=1}^m(1-\frac{1}{r_i})$. 
Thus $m=1$, which implies that $H_t$ is the only connected component of $\Supp M_t$ which intersects $G_t$. 

We show that $H_t\cup G_t$ is a $K_{\Ff_t}$-chain. Since $Q_t\cdot G_t<0$ and $Q_t$ is pseudo-effective, $G_t^2<0$. Since $G_t'$ is a smooth rational curve, $G_t$ is a smooth rational curve. Since $H_t$ is a string and $H_{l,t}\cdot G_t=1$, we deduce that $H_t\cup G_t$ is a string. Since $Z(\Ff_t',G_t')=1$ and $CS(\Ff_t',G_t')<0$, by Lemma \ref{lem: eigenvalue and cs and z formula} and Lemma \ref{lem: F chain and KF chain}(5), $Z(\Ff_t,G_t)=2$ and $G_t$ does not contain any Poincar\'{e}-Dulac singularity of $\Ff_t$. Therefore $H_t\cup G_t$ is an $\Ff_t$-chain. By our previous statements, we deduce that $H_t\cup G_t$ is a $(K_{\Ff_t}+\Delta_t)$-chain. Since $(H_t\cup G_t)\cap\Supp\Delta_t=0$, by applying Lemma \ref{lem: Intersection surface}, $H_t\cup G_t$ is a $K_{\Ff_t}$-chain.
\end{proof}

\noindent\textit{Proof of Lemma \ref{lem: deform zariski decomposition foliated pairs} continued}.
Let $U_2\subset U_1$ be an open neighborhood of $s$ and $\omega$ a hermitian metric on $X$, such that $\omega|_{X_t}$ is a K\"{a}hler metric for any $t\in U_2$. If $Q_{t}$ is nef for any $t\in U_2$, we finish the proof by letting $U=U_2$. Otherwise, there is a closed point $t_0\in U_2$ and an irreducible curve $G_{t_0}\subset X_{t_0}$ such that $Q_{t_0}\cdot G_{t_0}<0$. By Lemma \ref{lem: kf chain deforms} and Claim \ref{claim: gt is kft chain}, there is an open neighborhood $U_{t_0}\subset U_2$ of $t_0$ such that $G_{t_0}$ deforms over $U_{t_0}$. 

Suppose that $G_{t_0}$ deforms to $G_t$ for any $t\in U_{t_0}$. Then for any $t\in U_{t_0}$, $G_t\cdot [\omega|_{X_t}]$ is a constant. Thus we may extend $G_{t_0}$ to an effective divisor $G$ on $X$. Let $G_s:=G|_{X_s}$, then
$$0\leq P_s\cdot G_s=Q_s\cdot G_s=Q_{t_0}\cdot G_{t_0}<0,$$
a contradiction.
\end{proof}

\section{Proof of the main theorems}
\begin{thm}\label{thm: weak foliated surface pair inv plu} 
Assume that 
\begin{itemize}
    \item $X$ is a smooth threefold, 
    \item $\Ff$ is a rank 1 foliation on $X$,
    \item $\pi: X\rightarrow T$ is a smooth morphism to a curve,
    \item $s\in T$ is a closed point, 
    \item $\Delta$ is an effective $\mathbb Q$-Cartier $\mathbb Q$-divisor on $X$, and
    \item $A$ is an ample$/T$ $\Qq$-divisor on $X$ which does not contain any fiber of $\pi$
\end{itemize}
satisfying the following:
\begin{itemize}
    \item  $(X_t,\Ff_t)|_{t\in T}$ is a smooth family of canonical foliations of surfaces associated with $\pi$,
    \item $\Delta_t:=\Delta|_{X_t}$ and $A_t:=A|_{X_t}$ for any closed point $t\in T$,
 \item for any closed point $t\in T$ and any irreducible component $\Delta^0$ of $\Delta$, $\Delta^0|_{X_t}$ is irreducible and reduced, 
\item $K_{\Ff}+\Delta$ is $\pi$-pseudo-effective, and
\item $\Delta_s\wedge\textbf{\rm \textbf{B}}_{-}(K_{\Ff_s}+\Delta_s)=0$,
\end{itemize}
then there exists a rational number $0<\epsilon\ll 1$ and an open neighborhood $U$ of $s$, such that for any sufficiently divisible integer $m>0$,
\begin{enumerate}
\item $$h^0(X_t, m(K_{\Ff_t}+\Delta_t+\epsilon A_t))=\const,$$
for any $t\in U$, and
\item if 
\begin{itemize}
\item $K_{\Ff}+\Delta$ is big over $T$, and 
\item $\Delta_s\wedge\textbf{\rm \textbf{B}}_{+}(K_{\Ff_s}+\Delta_s)=0$,
\end{itemize}
then
$$h^0(X_t, m(K_{\Ff_t}+\Delta_t))=\const$$
for any $t\in U$.
\end{enumerate}
\end{thm}

\begin{proof}[Proof of Theorem \ref{thm: weak foliated surface pair inv plu}] By Lemma Lemma \ref{lem: contract -1curve in a chain}, Lemma \ref{lem: contract artificial chain step by step}, and Lemma \ref{lem: inv plur foliated pairs under bir contraction}, there is an open neighborhood $U_0$ of $s$ and a birational contraction 
$$\nu: (\pi^{-1}(U_0),\Ff|_{\pi^{-1}(U_0)},\Delta|_{\pi^{-1}(U_0)})\rightarrow (X',\Ff',\Delta')$$ 
over $U_0$, such that
\begin{itemize}
\item $A':=\nu_*A|_{\pi^{-1}(U_0)}$,
    \item any $(K_{\Ff'_s}+\Delta'_s)$-chain or any $(K_{\Ff'_s}+\Delta'_s+\epsilon A‘)$-chain does not contain any $(-1)$-curve for any $0<\epsilon\ll 1$, and
    \item for any $t\in U_0$, any rational number $0<\epsilon\ll 1$ and any sufficiently divisible integer $m>0$, $$h^0(X'_t, m(K_{\Ff'_t}+\Delta'_t+\epsilon A'_t))=h^0(X_t, m(K_{\Ff_t}+\Delta_t+\epsilon A_t))$$
and
$$h^0(X'_t, m(K_{\Ff'_t}+\Delta'_t))=h^0(X_t, m(K_{\Ff_t}+\Delta_t)).$$
\end{itemize}

Thus possibly replacing $T$ with $U_0$, $(X,\Ff,\Delta)$ with $(X',\Ff',\Delta')$ and $A$ with $A'$, we may assume that every $(K_{\Ff_s}+\Delta_s)$-chain and every $(K_{\Ff_s}+\Delta_s+\epsilon A_s)$-chain does not contain any $(-1)$-curve.

Fix a rational number $0<\epsilon\ll 1$. By Lemma \ref{lem: deform zariski decomposition foliated pairs}, there exists an open neighborhood $U_1$ of $s$ and four $\mathbb Q$-divisors $P^A,P,N^A$ and $N$ on $X$ over $U_1$, such that 
\begin{itemize}
\item $K_\Ff+\Delta+\epsilon A=P^A+N^A$,
\item $K_\Ff+\Delta=P+N$,
\item for any $t\in U_1$, $P^A_t:=P^A|_{X_t}$ and $N^A_t:=N^A|_{X_t}$ are the positive and the negative part of the Zariski decomposition of $K_{\Ff_t}+\Delta_t+\epsilon A_t$ respectively, 
\item for any $t\in U_1$, $P_t:=P|_{X_t}$ and $N_t:=N|_{X_t}$ are the positive and the negative part of the Zariski decomposition of $K_{\Ff_t}+\Delta_t$ respectively, and
\item intersection numbers deform over $U_1$.
\end{itemize} 

Let $m>0$ be a sufficiently divisible integer. By Lemma \ref{lem: kv vanishing foliation version 1}, for any integer $i>0$, $H^i(mP^A_s)=0$. Thus for any $t\in U_1$,
$$h^0(X_s, m(K_{\Ff_s}+\Delta_s+\epsilon A_s))=h^0(X_s, mP^A_s)=\chi(X_s, mP^A_s)=\chi(X_t, mP^A_t).$$
Since $P^A_t$ is nef and big, by Serre duality,
$$h^2(X_t, mP^A_t)=0,$$ 
which implies that
$$\chi(X_t, mP^A_t)\leq h^0(X_t, mP^A_t)=h^0(X_t, m(K_{\Ff_t}+\Delta_t+\epsilon A_t)).$$
Therefore
$$h^0(X_t, m(K_{\Ff_t}+\Delta_t+\epsilon A_t))\geq h^0(X_s, m(K_{\Ff_s}+\Delta_s+\epsilon A_s))$$
for any $t\in U_1$. Since $t\rightarrow h^0(X_t, m(K_{\Ff_t}+\Delta_t+\epsilon A_t))$ is an upper-semicontinuous function, there exists an open neighborhood $U_2\subset U_1$ of $s$, such that
$$h^0(X_t, m(K_{\Ff_t}+\Delta_t+\epsilon A_t))=\const$$
for any $t\in U_2$, which implies (1).

Under the assumptions of (2), by Lemma \ref{lem: lem: kv vanishing foliation version 2}, $H^i(mP_s)=0$ for any integer $i>0$. Thus for any $t\in U_1$,
$$h^0(X_s, m(K_{\Ff_s}+\Delta_s))=h^0(X_s, mP_s)=\chi(X_s, mP_s)=\chi(X_t, mP_t).$$
Since $K_{\Ff}+\Delta$ is big over $T$, $P_t$ is big for any $t\in U_1$. By Serre duality, 
$$h^2(X_t, mP_t)=0,$$
for any $t\in U_1$, which implies that
$$\chi(X_t, mP_t)\leq h^0(X_t, mP_t)=h^0(X_t, m(K_{\Ff_t}+\Delta_t)).$$
Since $t\rightarrow h^0(X_t, m(K_{\Ff_t}+\Delta_t))$ is an upper-semicontinuous function, there exists an open neighborhood $U_3\subset U_1$ of $s$, such that
$$h^0(X_t, m(K_{\Ff_t}+\Delta_t))=\const$$
for any $t\in U_3$, which implies (2).
\end{proof}

\begin{proof}[Proof of Theorem \ref{thm: foliated surface pair inv plur}] Let $X_s$ be a closed fiber of $\pi$ and $A$ a general ample$/T$ $\mathbb Q$-divisor. Since $K_{\Ff}+\Delta$ is big over $T$, we may pick a rational number $0<\epsilon\ll 1$ such that $K_\Ff+\Delta-\epsilon A$ is pseudo-effective over $T$. For any $t\in T$, suppose that $P_{t,\epsilon}$ and $N_{t,\epsilon}$ are the positive and the negative part of the Zariski decomposition of $K_{\Ff_t}+\Delta_t-\epsilon A_t$ respectively. We define
$$\Theta_s:=(\Delta_s-\epsilon A_s)-\Delta_s\wedge N_{s,\epsilon},$$
and let $\Theta$ be the unique $\Qq$-divisor on $X$ such that  $\Theta_s=\Theta|_{X_s}$. Since $\Theta_s+\epsilon A_s\geq 0$, $\Theta+A\geq 0$. Let $\Theta_t:=\Theta|_{X_t}$ for any $t\in T$. Then
$$K_{\Ff_s}+\Theta_s=P_{s,\epsilon}+(N_{s,\epsilon}-\Delta_s\wedge N_{s,\epsilon})$$ 
is pseudo-effective and
$$\textbf{B}_{+}(K_{\Ff_s}+\Theta_s+\epsilon A_s)\subset \textbf{B}_{-}(K_{\Ff_s}+\Theta_s).$$
Since
$$(\Delta_s-\Delta_s\wedge N_{s,\epsilon})\wedge(N_{s,\epsilon}-\Delta_s\wedge N_{s,\epsilon})=0,$$ 
we have
$$(\Theta_s+\epsilon A_s)\wedge\textbf{B}_{-}(K_{\Ff_s}+\Theta_s)=0,$$ 
which implies that
$$(\Theta_s+\epsilon A_s)\wedge\textbf{B}_{+}(K_{\Ff_s}+\Theta_s+\epsilon A_s)=0.$$
By Theorem \ref{thm: weak foliated surface pair inv plu}(2), there is an open neighborhood $U$ of $s$, such that
$$h^0(X_t,m(K_{\Ff_t}+\Theta_t+\epsilon A_t))=h^0(X_s,m(K_{\Ff_s}+\Theta_s+\epsilon A_s))$$
for any sufficiently divisible integer $m>0$ and any $t\in U$. 

Since $N_{s,\epsilon}\geq\Delta_s\wedge N_{s,\epsilon}\geq 0$ and $\Delta\geq\Theta+\epsilon A\geq 0$, we have 
\begin{align*}
|m(K_{\Ff_s}+\Delta_s)|&=|m((K_{\Ff_s}+\Delta_s-\epsilon A_s)+\epsilon A_s)|\\
&=|m((K_{\Ff_s}+\Delta_s-\epsilon A_s)-\Delta_s\wedge N_{s,\epsilon}+\epsilon A_s)|\\
&=|m(K_{\Ff_s}+\Theta_s+\epsilon A_s)|\\
&=|m(K_{\Ff}+\Theta+\epsilon A)|_{\pi^{-1}(U)}|_s\\
&\subset |m(K_\Ff+\Delta)|_{\pi^{-1}(U)}|_s\subset |m(K_{\Ff_s}+\Delta_s)|.
\end{align*}
Thus $|m(K_{\Ff_s}+\Delta_s)|=|m(K_\Ff+\Delta)|_{\pi^{-1}(U)}|_s$. In particular, for any $t\in U$,
$$h^0(X_t,m(K_{\Ff_t}+\Delta_t))=\const.$$
Since $X_s$ can be any closed fiber of $X\rightarrow T$, we have
$$h^0(X_t,m(K_{\Ff_t}+\Delta_t))=\const$$
for any $t\in T$.
\end{proof}

To prove Corollary \ref{cor: foliated surface pair vol inv}, we first state and prove another theorem on invariance of plurigenera.

\begin{thm}\label{thm: deformation invariance foliated pairs under scaling}  
Assume that 
\begin{itemize}
    \item $X$ is a smooth threefold, 
    \item $\Ff$ is a rank 1 foliation on $X$,
    \item $\pi: X\rightarrow T$ is a smooth morphism to a curve,
    \item $s\in T$ is a closed point, 
    \item $\Delta$ is an effective $\mathbb Q$-Cartier $\mathbb Q$-divisor on $X$, and
    \item $A$ is an ample$/T$ $\Qq$-divisor on $X$ which does not contain any fiber of $\pi$.
\end{itemize}
satisfying the following:
\begin{itemize}
    \item  $(X_t,\Ff_t)|_{t\in T}$ is a smooth family of canonical foliations of surfaces associated with $\pi$,
    \item $\Delta_t:=\Delta|_{X_t}$ and $A_t:=A|_{X_t}$ for any closed point $t\in T$,
 \item for any closed point $t\in T$ and any irreducible component $\Delta^0$ of $\Delta$, $\Delta^0|_{X_t}$ is irreducible and reduced, 
\item $K_{\Ff}+\Delta$ is $\pi$-pseudo-effective,
\end{itemize}
then there exists a rational number $0<\epsilon\ll 1$, such that for any sufficiently divisible integer $m>0$, 
$$h^0(X_t, m(K_{\Ff_t}+\Delta_t+\epsilon A_t))=\const.$$

\end{thm}

\begin{proof} Let $X_s$ be a closed fiber of $\pi$. For any closed point $t\in T$, let $P_t$ and $N_t$ be the positive part and the negative part of the Zariski decomposition of $K_{\Ff_t}+\Delta_t$ respectively. We define $\Theta_s:=\Delta_s-\Delta_s\wedge N_s$ and let $\Theta$ be the unique $\Qq$-divisor on $X$ such that $\Theta_s=\Theta|_{X_s}$. Let $\Theta_t:=\Theta|_{X_t}$ for any $t\in T$.

By Lemma \ref{lem: deform zariski decomposition foliated pairs}(1), there exists a rational number $0<\epsilon\ll 1$, an open neighborhood $U$ of $s$, two $\mathbb Q$-divisors $P^{A}$ and $N^{A}$ on $\pi^{-1}(U)$, such that 
\begin{itemize}
\item $(K_\Ff+\Theta+\epsilon A)|_{\pi^{-1}(U)}=P^{A}+N^{A}$, and
 \item for any $t\in U$, $P^{A}_{t}:=P^{A}|_{X_t}$ and $N^{A}_t:=N^{A}|_{X_t}$ are the positive and the negative part of the Zariski decomposition of $K_{\Ff_t}+\Theta_t+\epsilon A_t$ respectively.
\end{itemize}

For any sufficiently divisible integer $m>0$, we have
$$|m(K_{\Ff_s}+\Delta_s+\epsilon A_s)|=|m(K_{\Ff_s}+\Theta_s+\epsilon A_s)|=|mP^{A}_s|.$$
By Theorem \ref{thm: weak foliated surface pair inv plu}(1), there exists an open neighborhood $U_1\subset U$ of $s$, such that $h^0(X_s, m(K_{\Ff_s}+\Theta_s+\epsilon A_s))=h^0(X_t, m(K_{\Ff_t}+\Theta_t+\epsilon A_t))$ for any $t\in U_1$. Thus for any sufficiently divisible integer $m>0$, 
\begin{align*}
|m(K_{\Ff_s}+\Delta_s+\epsilon A_s)|&=|m(K_{\Ff_s}+\Theta_s+\epsilon A_s)|\\
&=|m(K_{\Ff}+\Theta+\epsilon A)|_{\pi^{-1}(U_1)}|_s\\
&\subset |m(K_{\Ff}+\Delta+\epsilon A)|_{\pi^{-1}(U_1)}|_s\subset |m(K_\Ff+\Delta+\epsilon A)|_s,
\end{align*}
which implies that
$$|m(K_{\Ff_s}+\Delta_s+\epsilon A_s)|\subset |m(K_\Ff+\Delta+\epsilon A)|_s.$$
Thus
$$|m(K_{\Ff_s}+\Delta_s+\epsilon A_s)|=|m(K_\Ff+\Delta+\epsilon A)|_s.$$
Thus for any $t\in U_1$, 
$$h^0(X_t, m(K_{\Ff_t}+\Delta_t+\epsilon A_t))=\const.$$ 
Since $X_s$ can be any closed fiber of $X\rightarrow T$, we have
$$h^0(X_t, m(K_{\Ff_t}+\Delta_t+\epsilon A_t))=\const$$
for any $t\in T$. 
\end{proof}

\begin{proof}[Proof of Corollary \ref{cor: foliated surface pair vol inv}] First we deal with the case when $K_\Ff+\Delta$ is $\pi$-pseudo-effective. Let $A$ be a general ample$/T$ $\mathbb Q$-divisor on $X$ and $X_s$ a closed fiber of $\pi$. By Theorem \ref{thm: deformation invariance foliated pairs under scaling}, for any $t\in T$,
\begin{align*}
\vol(X_t,K_{\Ff_t}+\Delta_t)&=\lim_{\epsilon\rightarrow+0}\vol(X_t,K_{\Ff_t}+\Delta_t+\epsilon A_t)\\
&=\lim_{\epsilon\rightarrow+0}\vol(X_s,K_{\Ff_s}+\Delta_s+\epsilon A_s)=\vol(X_s,K_{\Ff_s}+\Delta_s),
\end{align*}
and the proof is finished.

Now we deal with the case when $K_\Ff+\Delta$ is not $\pi$-pseudo-effective. We may suppose that there exists $s\in T$ such that $K_{\Ff_s}+\Delta_s$ is pseudo-effective, otherwise $\vol(X_t,K_{\Ff_t}+\Delta_t)=0$ for any $t\in T$ and there is nothing to prove.

Let $A$ be a general ample$/T$ $\mathbb Q$-divisor. Let $\lambda$ be the $\pi$-pseudo-effective threshold of $K_\Ff+\Delta$ with respect of $A$, i.e. 
$$\lambda:=\inf\{r\in\mathbb R|K_\Ff+\Delta+rA \text{ is }\pi\text{-pseudo-effective}\}.$$
Since $K_\Ff+\Delta$ is not $\pi$-pseudo-effective, $\lambda>0$. 

For any integer $i>0$, let $\lambda_i:=\lambda+\frac{1}{i}$. Since $K_{\Ff_s}+\Delta_s$ is pseudo-effective, $K_{\Ff_s}+\Delta_s+\lambda A_s$ is big. Thus 
$$\vol(K_{\Ff_s}+\Delta_s+\lambda_i A_s)\geq\vol(K_{\Ff_s}+\Delta_s+\lambda A_s)>0.$$ 
However, by what we have already proved, 
\begin{align*}
0&<\vol(K_{\Ff_s}+\Delta_s+\lambda A_s)=\lim_{i\rightarrow+\infty}\vol(K_{\Ff_s}+\Delta_s+\lambda_i A_s)\\
&=\lim_{i\rightarrow+\infty}\vol(K_{\Ff_t}+\Delta_t+\lambda_i A_t)=\vol(K_{\Ff_t}+\Delta_t+\lambda A_t)
\end{align*}
for any $t\in T$. 

Since $\vol(K_{\Ff_t}+\Delta_t+rA_t)$ is a continuous function for $(r,t)\in\mathbb R\times T$, we may pick $\lambda'<\lambda$ such that $\vol(K_{\Ff_t}+\Delta_t+\lambda' A_t)>0$ for any $t\in T$. Thus $\lambda$ is not the $\pi$-pseudo-effective threshold of $K_\Ff+\Delta$ with respect to $A$, a contradiction.\end{proof}


\begin{thebibliography}{99}

\bibitem[AD13]{AD13} C.Araujo and S.Druel, \textit{On Fano foliations}, Adv. Math. 238 (2013), 70-118.

\bibitem[Bru00]{Bru00} M.Brunella, \textit{Birational geometry of foliations}, Monograf\'{i}as de Matem\'atica, Instituto de Matem\'atica Pura e Aplicada(IMPA), Rio de janeiro, 2000.

\bibitem[Bru01]{Bru01} M.Brunella, \textit{Invariance par d\'eformations de la dimension de Kodaira d'un feuilletage sur une surface}, Essays on geometry and related topics, Vol.1,2, Monogr. Enseign. Math., vol.38, Enseignement Math., Geneva, 2001, pp.113-132.

\bibitem[Bru02]{Bru02} M.Brunella, \textit{Foliations on complex projective surfaces}, arXiv:math/0212082, 2002, 31p.

\bibitem[CF15]{CF15} P.Cascini and E.Floris, \textit{On invariance of plurigenera for foliations on surfaces}, to appear in Crelle's Journal, arXiv:math/1502.00817, 2015, 38p.

\bibitem[ELMNP06]{ELMNP06} L.Ein, R.Lazarsfeld, M.Musta\c{t}\v{a}, M.Nakayame, and M.Popa, \textit{Asymptotic invariants of base loci}, Ann. Inst. Fourier, \textbf{56}, 2006, pp.1701-1734.

\bibitem[ELMNP09]{ELMNP09} L.Ein, R.Lazarsfeld, M.Musta\c{t}\v{a}, M.Nakayame, and M.Popa, \textit{Restricted volumes and asymptotic intersection theory}, Amer. J. Math, \textbf{131}, 2009, no. 3, pp. 607-651.

\bibitem[FM94]{FM94} R.Friedman and J.W.Morgan, \textit{Smooth four-manifolds and complex surfaces}, Ergebnisse der Mathematik und ihrer Grenzgebiete (3) [Results in Mathematics and Related Areas (3)], vol. 27, Springer-Verlag, Berlin, 1994.

\bibitem[HMX13]{HMX13} C.D.Hacon, J.M$^c$kernan and C.Xu, \textit{On the birational automorphisms of varieties of general type}, Ann. of Math. 177 (2013), no. 3, pp. 1077-1111.

\bibitem[HMX14]{HMX14} C.D.Hacon, J.M$^c$kernan and C.Xu, \textit{Boundness of moduli of varieties of general type}, arXiv preprint arXiv:1412.1186, 2014, 37p.

\bibitem[KM98]{KM98} J.Koll\'{a}r and S.Mori, \textit{Birational geometry of algebraic varieties}, Cambridge tracts in mathematics, vol. 134, Cambridge University Press, 1998.

\bibitem[Kod63]{Kod63} K.Kodaira, \textit{On stability of compact submanifolds of complex manifolds}, Amer. J. Math. 85, pp. 79-94, 1963.

\bibitem[Laz04]{Laz04} R.Lazarsfeld, \textit{Positivity in algebraic geometry I,II}, Ergebnisse der Mathematik und ihrer Grenzgebiete. 3. Folge. A Series of Modern Surveys in Mathematics [Results in Mathematics and Related Areas. 3rd Series. A Series of Modern Surveys in Mathematics], vol. 48, ``Classical Setting: Line Bundles and Linear Series" and vol.49, ``Positivity for vector bundles, and multiplier ideals", Springer-Verlag, Berlin, 2004.

\bibitem[LPT11]{LPT11} F.Loray, J.Pereira, and F.Touzet, \textit{Singular foliations with trivial canonical class}, arXiv eprint arXiv:1107.1538, 2011, 39p.

\bibitem[McQ08]{McQ08} M.McQuillan, \textit{Canonical models of foliations}. Pure Appl. Math. Q., 4(3, part 2), 2008, pp. 877-1012.

\bibitem[Mum61]{Mum61} D.Mumford, \textit{The topology of normal singularities on algebraic surfaces}, Publ. Math. I.H.E.S. \textbf{9}, 1961, pp. 5-22.

\bibitem[Nak00]{Nak00} M.Nakamaye, \textit{Stable base loci of linear series}, Math. Ann. \textbf{318}, 2000, no.4, 837-847.

\bibitem[Siu98]{Siu98} Y.-T.Siu, \textit{Invariance of plurigenera}, Invent. Math., 134(3), 1998, pp. 661-673

\bibitem[Spi17]{Spi17} Calum Spicer, \textit{Higher dimensional foliated mori theory}, arXiv:1709.06850v3, 2017.

\end{thebibliography}
\end{document}